\documentclass[a4paper,11pt,reqno,oneside]{amsart}
\usepackage{amsmath,amssymb,amsthm,mathtools}
\usepackage{mathrsfs,enumitem,xcolor}
\usepackage[scale=0.75,hcentering]{geometry}
\usepackage{float}
\usepackage{microtype,aliascnt,hyperref}
\usepackage{array,booktabs,longtable,multirow}
\usepackage[capitalize,nameinlink,noabbrev]{cleveref}
\makeatletter
\@ifundefined{subjclassname@2020}{\@namedef{subjclassname@2020}{%
\textup{2020} Mathematics Subject Classification}}{}
\makeatother
\definecolor{cobalt}{rgb}{0.0,0.28,0.67}
\hypersetup{bookmarksopen=true,
pdftitle={Upper Bounds on the Strassen CP Rank Gap and Symmetric Rank Gap of Third-Order Tensors},
pdfauthor={Jianze Li},
pdfsubject={CP and symmetric ranks over the real and complex fields},
pdfkeywords={symmetric tensor, CP rank, Waring rank, Comon rank gap, direct sum},
colorlinks=true,linkcolor=cobalt,citecolor=cobalt,urlcolor=cobalt}
\theoremstyle{plain}
\newtheorem{theorem}{Theorem}[section]
\newaliascnt{corollary}{theorem}
\newtheorem{corollary}[corollary]{Corollary}
\aliascntresetthe{corollary}
\newaliascnt{lemma}{theorem}
\newtheorem{lemma}[lemma]{Lemma}
\aliascntresetthe{lemma}
\newaliascnt{proposition}{theorem}
\newtheorem{proposition}[proposition]{Proposition}
\aliascntresetthe{proposition}
\theoremstyle{definition}
\newaliascnt{definition}{theorem}

\aliascntresetthe{definition}
\newaliascnt{problem}{theorem}
\newtheorem{problem}[problem]{Open Problem}
\aliascntresetthe{problem}
\theoremstyle{remark}
\newaliascnt{remark}{theorem}
\newtheorem{remark}[remark]{Remark}
\aliascntresetthe{remark}
\newaliascnt{example}{theorem}

\aliascntresetthe{example}
\crefname{theorem}{Theorem}{Theorems}
\Crefname{theorem}{Theorem}{Theorems}
\crefname{corollary}{Corollary}{Corollaries}
\Crefname{corollary}{Corollary}{Corollaries}
\crefname{lemma}{Lemma}{Lemmas}
\Crefname{lemma}{Lemma}{Lemmas}
\crefname{problem}{Open Problem}{Open Problems}
\Crefname{problem}{Open Problem}{Open Problems}
\crefname{proposition}{Proposition}{Propositions}
\Crefname{proposition}{Proposition}{Propositions}
\crefname{definition}{Definition}{Definitions}
\Crefname{definition}{Definition}{Definitions}
\crefname{remark}{Remark}{Remarks}
\Crefname{remark}{Remark}{Remarks}
\crefname{example}{Example}{Examples}
\Crefname{example}{Example}{Examples}
\newcommand{\tens}[1]{\boldsymbol{#1}}
\newcommand{\matr}[1]{\mathbf{#1}}
\newcommand{\vect}[1]{\boldsymbol{#1}}
\newcommand{\set}[1]{\mathcal{#1}}
\newcommand{\eqdef}{\triangleq}
\newcommand{\RR}{\mathbb{R}}
\newcommand{\CC}{\mathbb{C}}

\newcommand{\FF}{\mathbb{F}}
\newcommand{\Sym}{\mathsf{Sym}}
\newcommand{\Span}{\mathsf{span}}
\newcommand{\rank}{\mathsf{rank}}
\newcommand{\Ess}{\mathsf{Ess}}

\newcommand{\supp}{\mathsf{supp}}
\newcommand{\id}{\mathsf{id}}
\numberwithin{equation}{section}
\allowdisplaybreaks
\setlist[enumerate]{label=\textup{(\roman*)},leftmargin=*}
\title[Upper Bounds on the Strassen CP Rank Gap and Symmetric Rank Gap]
{Upper Bounds on the Strassen CP Rank Gap and Symmetric Rank Gap of Third-Order Tensors}
\author{Jianze Li}
\address{School of Science, Sun Yat-sen University, Shenzhen, Guangdong, China}
\email{lijianze@mail.sysu.edu.cn}
\date{\today}

\begin{document}
\raggedbottom

\begin{abstract}
For two third-order tensors over $\RR$ or $\CC$, we define the Strassen CP rank gap as the sum of their CP ranks minus the
CP rank of their direct sum.
For arbitrary third-order tensors, we derive explicit upper bounds
on this gap in terms of CP rank excess and dimension.
For symmetric third-order tensors, we similarly define the Strassen symmetric rank gap using symmetric rank and establish
explicit upper bounds in terms of symmetric rank excess, CP rank
excess, and dimension.
These bounds also yield sufficient conditions for the respective gaps
to vanish, ensuring additivity of CP rank or symmetric rank under
direct sums.
\end{abstract}
\maketitle
\enlargethispage{6pt}

\section{Introduction}\label{sec:introduction}

Tensor CP rank is the minimum number of simple tensors needed to represent
a tensor \cite{hitchcock1927expression,kolda2009tensor,kruskal1977three,landsberg2012tensors}.   For symmetric tensors,
one may further require all simple summands to be symmetric \cite{comon2008symmetric}, and the
resulting symmetric rank agrees with the Waring rank of the associated
homogeneous polynomial
\cite{comon2008symmetric,landsberg2012tensors}. 
The notion goes back to Hitchcock
\cite{hitchcock1927expression} and now plays a central
role in multilinear algebra, algebraic complexity, signal processing,
statistics, and data analysis; see, for example,
\cite{comon1994independent,Como10:book,kolda2009tensor,landsberg2012tensors,LUC2017globally,ULC2019}. 

In this paper, the base field is $\FF\in\{\RR,\CC\}$ unless
stated otherwise, and every vector space is finite-dimensional over
$\FF$.  Let $\set{U}_1,\set{U}_2,\set{U}_3$ be such vector spaces and let
$\tens{T}\in\set{U}_1\otimes\set{U}_2\otimes\set{U}_3$.  Its CP rank is
\[
r_\FF(\tens{T})\eqdef
\min\left\{q:\
\tens{T}=\sum_{h=1}^{q}
\vect{a}_h\otimes\vect{b}_h\otimes\vect{c}_h\right\},
\]
where $\vect{a}_h\in\set{U}_1$, $\vect{b}_h\in\set{U}_2$, and
$\vect{c}_h\in\set{U}_3$.  To record the intrinsic size of each mode,
fix $\{i,j,k\}=\{1,2,3\}$ and contract in modes $j$ and $k$.  This
gives the mode-$i$ flattening
$\mathsf{Flat}_i(\tens{T}):
\set{U}_j^*\otimes\set{U}_k^*\longrightarrow\set{U}_i.$ 
Its image is the mode-$i$ support of $\tens{T}$.  Accordingly, we
define
$d_i(\tens{T})\eqdef\rank\bigl(\mathsf{Flat}_i(\tens{T})\bigr)$
and $d(\tens{T})\eqdef\max_{1\le i\le3}d_i(\tens{T}).$
We call $d_i(\tens{T})$ the mode-$i$ rank and $d(\tens{T})$ the maximum mode rank. 
It holds that 
$d(\tens{T})\le r_\FF(\tens{T}).$  
Let
$\tens{T}_1\in\set{U}_1\otimes\set{U}_2\otimes\set{U}_3$ and $\tens{T}_2\in\set{V}_1\otimes\set{V}_2\otimes\set{V}_3$,  
where $\set{V}_1,\set{V}_2,\set{V}_3$ are also finite-dimensional vector
spaces over $\FF$.  
After placing the
factor spaces of $\tens{T}_1$ and $\tens{T}_2$ in direct sum, their 
direct sum tensor
$$\tens{T}_1\oplus\tens{T}_2\in
(\set{U}_1\oplus\set{V}_1)\otimes
(\set{U}_2\oplus\set{V}_2)\otimes
(\set{U}_3\oplus\set{V}_3)$$ 
has $\tens{T}_1$ and $\tens{T}_2$ as its two diagonal blocks, while all mixed
blocks vanish.  Therefore, concatenating minimum decompositions gives
$r_\FF(\tens{T}_1\oplus\tens{T}_2)
\le r_\FF(\tens{T}_1)+r_\FF(\tens{T}_2).$  
In the symmetric setting, the three factor spaces are all equal to a vector space $\set{U}$.
For a symmetric tensor 
$\tens{S}\in\Sym^3(\set{U})$, its symmetric rank is
\[
s_\FF(\tens{S})\eqdef
\min\left\{q:\
\tens{S}=\sum_{h=1}^{q}\lambda_h\vect{u}_h^{\otimes3},\ 
\lambda_h\in\FF,\ \vect{u}_h\in\set{U}\right\}.
\]
In this case, the three actual mode spaces coincide.  Their common
value is the essential space $\Ess(\tens{S})$; hence
$d(\tens{S})=\dim\Ess(\tens{S})$, which is called the essential dimension of $\tens{S}$.  
The ranks satisfy
$d(\tens{S})\le r_\FF(\tens{S})\le s_\FF(\tens{S}).$ 

The mode ranks provide lower bounds for CP rank; in the symmetric case,
their common value also bounds symmetric rank.  We thus define the CP rank excess, the
symmetric rank excess, and the Comon rank gap by
\[
e_{r,\FF}(\tens{T})\eqdef r_\FF(\tens{T})-d(\tens{T}),
\ \ 
e_{s,\FF}(\tens{S})\eqdef s_\FF(\tens{S})-d(\tens{S}),
\ \ 
\Delta_{c,\FF}(\tens{S})\eqdef s_\FF(\tens{S})-r_\FF(\tens{S}), 
\]
where the latter two quantities are defined only for symmetric tensors \cite{li2026comonrankgap,li2026hankeladditivity}. 
We also measure the possible failure of rank
additivity by the Strassen CP rank gap and Strassen symmetric rank gap \cite{li2026hankeladditivity}, defined as 
\begin{align*}
\Delta_{r,\FF}(\tens{T}_1,\tens{T}_2)
&\eqdef r_\FF(\tens{T}_1)+r_\FF(\tens{T}_2)
-r_\FF(\tens{T}_1\oplus\tens{T}_2), \\
\Delta_{s,\FF}(\tens{S}_1,\tens{S}_2)
&\eqdef s_\FF(\tens{S}_1)+s_\FF(\tens{S}_2)
-s_\FF(\tens{S}_1\oplus\tens{S}_2),
\end{align*}
where the second expression concerns symmetric tensors. 
All the above quantities are nonnegative.  For symmetric tensors,
they also satisfy
\begin{equation*}
\Delta_{s,\FF}(\tens{S}_1,\tens{S}_2) - \Delta_{r,\FF}(\tens{S}_1,\tens{S}_2)
= \Delta_{c,\FF}(\tens{S}_1)+\Delta_{c,\FF}(\tens{S}_2)
-\Delta_{c,\FF}(\tens{S}_1\oplus\tens{S}_2).
\end{equation*}
% For reference, \cref{tab:basic-notation} summarizes the notation used
% in the quantitative statements.  Every entry in the general-tensor
% column also applies to a symmetric tensor; the last column records
% only the additional or specialized symmetric notion. 
Since
the base field remains fixed within each statement and proof, we
henceforth suppress its subscript and write
$r,s,e_r,e_s,\Delta_c,\Delta_r$ and $\Delta_s$ for simplicity. 

% \begin{table}[H]
% \centering
% \small
% \setlength{\tabcolsep}{5pt}
% \renewcommand{\arraystretch}{1.2}
% \begin{tabular}{|>{\centering\arraybackslash}m{.15\linewidth}|
% >{\centering\arraybackslash}m{.36\linewidth}|
% >{\centering\arraybackslash}m{.39\linewidth}|}
% \hline
% Quantity
% & ordinary tensor $\tens{T}\in\set{U}_1\otimes\set{U}_2\otimes\set{U}_3$
% & Symmetric tensor $\tens{S}\in\Sym^3(\set{U})$ \\
% \hline
% \multirow{2}{*}{Rank}
% & \multicolumn{2}{|c|}{CP rank $r_\FF(\tens{T})$} \\
% \cline{2-3}
% & \multicolumn{1}{c|}{$\diagup$}
% & Symmetric rank $s_\FF(\tens{S})$ \\
% \hline
% Dimension
% & $d_i(\tens{T})$ and
% $d(\tens{T})$
% & $d(\tens{S})$ \\
% \hline
% \multirow{2}{*}{Rank excess}
% & \multicolumn{2}{|c|}{$e_{r,\FF}(\tens{T})$} \\
% \cline{2-3}
% & \multicolumn{1}{c|}{$\diagup$}
% & $e_{s,\FF}(\tens{S})$ \\
% \hline
% \multirow{2}{*}{Strassen CP rank gap}
% & \multicolumn{2}{|c|}{$\Delta_{r,\FF}(\tens{T}_1,\tens{T}_2)$} \\
% \cline{2-3}
% & \multicolumn{1}{c|}{$\diagup$}
% & $\Delta_{s,\FF}(\tens{S}_1,\tens{S}_2)$ \\
% \hline
% Comon rank gap
% & \multicolumn{1}{c|}{$\diagup$}
% & $\Delta_{c,\FF}(\tens{S})$ \\
% \hline
% \end{tabular}
% \caption{Rank, dimension, excess, and rank gap notations.}
% \label{tab:basic-notation}
% \end{table}

We now recall Strassen's direct sum conjecture.  In 1973, Strassen
asked whether CP rank is always additive under direct sums
\cite{strassen1973division}.  Feig and
Winograd subsequently studied the stronger question of whether every
minimum algorithm must separate into minimum algorithms for the two
blocks \cite{feig1984directsum}.  Several important
positive cases were then proved.  Ja'Ja' and Takche proved additivity
when one mode rank of either summand is at most two
\cite{jaja1986validity}; see also
\cite[Theorem~3.0.6]{rupniewski2024additivity}.  Buczy\'nski,
Postinghel, and Rupniewski obtained further criteria for small formats
and ranks \cite[Theorem~1.1]{buczynski2020additivity}. Rupniewski later proved, among other results, additivity over $\CC$
when both ranks are at most seven and a field-independent additivity criterion for tensors one of whose slice spaces has hook shape $(1,2)$
\cite[Theorem~1.0.3]{rupniewski2024additivity}. 
More recently, CP-rank additivity was established in
\cite{li2026hankeladditivity} for the direct sum of any complex
Hankel tensor of any order $m\ge3$ and any complex tensor of the
same order. 
For unrestricted third-order tensors, however, Shitov proved the existence of counterexamples over every infinite field
\cite[Theorem~4]{shitov2019directsum}.  
The quantitative estimates in
his construction imply that the gap is unbounded
\cite[Proposition~15, Lemma~16, and the proof of Theorem~7]
{shitov2019directsum}.  
In other words, over each field considered here, we have
$\sup_{\tens{T}_1,\tens{T}_2}\Delta_{r,\FF}(\tens{T}_1,\tens{T}_2)=\infty$ for 
$\FF\in\{\RR,\CC\}$,
where the supremum ranges over all third-order tensors.  
We refer to
\cite{casarotti2020comon} for a survey of earlier developments on Strassen's direct sum conjecture.  

As a symmetric analogue, Strassen's symmetric direct sum conjecture was explicitly formulated by Carlini, Catalisano, and
Chiantini
\cite{carlini2015symmetricstrassen}.
Over $\CC$, an earlier positive result established additivity for any number of monomials of the same degree at least two
in pairwise disjoint variable sets 
\cite[Proposition~3.1 and Theorem~3.2]{carlini2012monomials}.
Further cases include a pure-power summand, pairs of binary forms, and certain forms of small rank \cite[Proposition~3.1, Corollary~3.2, and Theorem~5.6]
{carlini2015symmetricstrassen}.  The framework of $e$-computability
then gave additivity for several infinite families and for any number
of forms $F_i$ of the same degree in disjoint variable sets that are
$e$-computable for the same $e$ and satisfy $(F_i^\perp)_e=0$
\cite[Definition~3.5 and Theorems~5.5 and~6.1]
{carlini2018ecomputability}, where $F_i^\perp$ denotes the apolar ideal of $F_i$ in the dual polynomial ring of its own variable set.  
Teitler obtained general criteria from
sharp catalecticant, affine-subspace, and singularity bounds; in the catalecticant case with derivatives of degree at least two, every minimum
Waring decomposition also separates
into decompositions of the summands
\cite[Theorems~1.9, 2.9, and~2.11]{teitler2015sufficient}.  Methods using
Waring loci yielded further strong separation results
\cite[Theorems~4.6 and~5.1]{carlini2017waringloci}. 
More recently, Shitov's preprint on higher rank substitutions reported
counterexamples to symmetric rank additivity under direct sums over every infinite field of
characteristic different from two and three
\cite[Theorem~3.9]{shitov2023higher}. 

These counterexamples in \cite{shitov2019directsum,shitov2023higher} show that the respective Strassen rank gaps need not vanish. 
This naturally leads to the following quantitative question:
\begin{quote}
\itshape
How large can the Strassen CP rank gap and symmetric rank gap be in terms of the rank excesses, ranks, and mode ranks of the summands?
\end{quote}
We study this question for ordinary and symmetric third-order
tensors over both $\RR$ and $\CC$, obtaining quantitative bounds
and sufficient conditions for additivity. 
For ordinary tensors, our main result bounds the Strassen CP rank gap
by an explicit function of the smaller CP rank excess.
This function grows at most linearly and vanishes throughout a
low-excess range. The bound is independent of the rank and format
of the other summand, so a sufficiently small excess on either side
guarantees additivity with every partner.
Combining this result with known maximum rank estimates gives bounds depending only on the smaller CP rank and quadratic bounds in the smaller
maximum mode rank. 
For symmetric tensors, we first establish a linear bound in the
smaller symmetric rank excess and refine it using both excesses
and the Comon rank gap of the other summand.
These estimates give additivity conditions that combine a small
excess on one side with equality, or a small difference, between
CP rank and symmetric rank on the other.
They also yield additivity criteria based only on symmetric ranks.
A further main result bounds the Strassen symmetric rank gap by a
piecewise linear function of the smaller CP rank excess. 
We also obtain quadratic bounds for the Strassen symmetric rank gap in terms
of the smaller essential dimension.
% The key input is a maximum-rank estimate for cubic forms that
% uses the generic complex rank in one fewer variable.
When both summands have at most four essential variables,
we prove that symmetric rank is additive.
% The precise bounds and additivity conditions are stated in
% \cref{sec:main-results}. 

The paper is organized as follows. \Cref{sec:preliminaries} collects the notation and
auxiliary results used in the proofs.
\Cref{sec:main-results} states the main gap bounds and
compares their additivity consequences with existing results. 
\Cref{sec:ordinary-proofs} proves
the CP rank estimates.  \Cref{sec:polar-tools,sec:rank-inputs,sec:gap-proofs}
develop the common tools and prove the symmetric excess
and dimension estimates.  \Cref{sec:conclusion} summarizes the
consequences and open directions. 

\section{Preliminaries}\label{sec:preliminaries}

\subsection{Notation and polarization}

We write ordinary tensors as $\tens{T}$ and symmetric tensors as $\tens{S}$,
matrices as $\matr{A}$, vectors as $\vect{a}$, and sets and vector spaces as
$\set{U}$.  Ordinary italic letters such as $F$ and $Q$ denote
polynomial forms.  All transposes $\matr{A}^{\mathsf T}$ are algebraic
rather than conjugate transposes.  The
operators $\rank$, $\Span$, $\Ess$, and $\supp$ denote matrix or map
rank, linear span, essential space, and support.  Our zero-tensor and
positive-part conventions are
$r(0)\eqdef0, s(0)\eqdef0, d(0)\eqdef0$
and 
$x_+\eqdef\max\{x,0\}.$ 
Let $\set{V}$ be a finite-dimensional vector space over $\FF$.  We
identify $\Sym^3(\set{V})$ with homogeneous cubic polynomials on $\set{V}^*$.
We write normalized symmetric products without a multiplication
symbol and use $\otimes$ for ordinary tensor products. 
For $F\in\Sym^3(\set{V})$ and $\alpha\in\set{V}^*$, put
$C_F(\alpha)\eqdef\alpha\mathbin{\lrcorner}F=\frac13D_\alpha F$,
where $D_\alpha$ denotes the directional derivative along $\alpha$ on $\set{V}^*$. 
Thus $C_F(\alpha)(\alpha,\alpha)=F(\alpha)$ and
$\alpha\mathbin{\lrcorner}\vect{u}^3
=\alpha(\vect{u})\vect{u}^2$.  The same notation $C_F(\alpha)$ denotes
the corresponding symmetric matrix in coordinates.  For any homogeneous
form $G$, we write $\Ess(G)$ for its smallest variable space.
Since every scalar in either field has a cube root, coefficients in cubic
Waring decompositions may be absorbed into the linear forms.
We record the following elementary properties of minimum decompositions
and omit their proofs; for the support statement in the CP case,
see \cite[Lemma~2.8]{buczynski2020additivity}.

\begin{lemma}[Minimum decompositions]\label{lem:minimum}
Let
$\tens{T}
=\sum_{j=1}^{r(\tens{T})}
\vect{a}_{j,1}\otimes\vect{a}_{j,2}\otimes\vect{a}_{j,3}$  
be a minimum CP decomposition. For each
$\set{I}\subseteq\{1,\ldots,r(\tens{T})\}$, define
$\tens{T}_{\set{I}}\eqdef
\sum_{j\in \set{I}}
\vect{a}_{j,1}\otimes\vect{a}_{j,2}\otimes\vect{a}_{j,3}.$ 
Then
\[
r(\tens{T}_{\set{I}})=|\set{I}|,
\ \ 
\mathsf{im}\bigl(\mathsf{Flat}_k(\tens{T}_{\set{I}})\bigr)
=\Span\{\vect{a}_{j,k}:j\in \set{I}\},
\ \  k=1,2,3.
\]
Moreover, for every pair $1\le k<\ell\le3$, we have 
$\dim\Span\{
\vect{a}_{j,k}\otimes\vect{a}_{j,\ell}:j\in \set{I}
\}=|\set{I}|.$ 
Similarly, let
$\tens{S}=\sum_{j=1}^{s(\tens{S})}\lambda_j\vect{u}_j^{\otimes3}$ with $\lambda_j\in\FF\setminus\{0\}$  
be a minimum Waring decomposition. For each
$\set{J}\subseteq\{1,\ldots,s(\tens{S})\}$, define
$\tens{S}_{\set{J}}\eqdef
\sum_{j\in \set{J}}\lambda_j\vect{u}_j^{\otimes3}.$ 
Then
\[
s(\tens{S}_{\set{J}})=|\set{J}|,
\ \ 
\Ess(\tens{S}_{\set{J}})=\Span\{\vect{u}_j:j\in \set{J}\}.
\]
\end{lemma}

\subsection{Rank bounds}

Set $g_3(0)\eqdef0$.  The Alexander--Hirschowitz theorem, in the form stated
in \cite[Theorem~2]{bernardi2018hitchhiker}, gives the generic
complex Waring rank of cubic forms in $n$ variables
\begin{equation}\label{eq:generic}
g_3(n)=
\begin{cases}
\lceil(n+1)(n+2)/6\rceil,&n\ge1,\ n\ne5,\\
8,&n=5.
\end{cases}
\end{equation}
The quantity $g_3(n)$ always refers to the complex generic Waring rank,
even in statements over $\RR$. Define the maximum symmetric rank by
$M_{3,\FF}(n)\eqdef
\max\{s(F):F\in\Sym^3(\FF^n)\}$ with 
$M_{3,\FF}(0)\eqdef0.$ 
Over both fields, we have
$M_{3,\FF}(n)\le2g_3(n)$ by
\cite[Theorems~1--3]
{blekherman2015maximum}.  We shall also use
\begin{equation}\label{eq:rank-cap}
M_{3,\FF}(n)\le2n-1,
\ \  1\le n\le4.
\end{equation}
For $n\le3$, this follows from the real and complex classifications in
\cite[Theorems~2.1--2.2]{banchi2015ternary}. 
For $n=4$ over $\CC$, see
\cite[Section~2]{seigal2020ranks}; over $\RR$, see
\cite[Proposition~5.2]{li2026comonrankgap}.

\subsection{Matrix spaces and tensor rank}

Let $\set{B}$ and $\set{C}$ be finite-dimensional vector spaces.  For a linear subspace
$\set{W}\subseteq\set{B}\otimes\set{C}$, define the rank of this matrix space by
\[
r(\set{W})\eqdef
\min\bigl\{q:\set{W}\subseteq\Span\{\matr{R}_1,\ldots,\matr{R}_q\},\ 
\rank(\matr{R}_j)=1\text{ for all }j\bigr\}.
\]
A rank-one cover of $\set{W}$ is a finite family of rank-one matrices whose
span contains $\set{W}$; it is minimum if its cardinality is $r(\set{W})$.
If $\tens{T}\in\set{A}\otimes\set{B}\otimes\set{C}$, then
$r(\tens{T})=r(\tens{T}(\set{A}^*))$, where
$\tens{T}(\set{A}^*)$ is the image of contraction in the first mode.
This is the slice-space description of tensor rank in
\cite[Lemma~2.9]{buczynski2020additivity}.  The members of a minimum rank-one
cover are linearly independent. 
We use the following support estimate \cite[Proposition~2.10(ii)]{buczynski2020additivity}.

\begin{lemma}[Ambient-support inequality]\label{lem:ambient-support}
Let $\set{W}\subseteq\set{B}_0\otimes\set{C}$ have
first-factor support $\set{B}_0\subseteq\set{B}$.
If a subspace $\set{R}\subseteq\set{B}\otimes\set{C}$
is spanned by rank-one matrices, contains $\set{W}$,
and has first-factor support $\set{B}_{\set{R}}$, then 
\[
\dim\set{R}\ge r(\set{W})+\dim\set{B}_{\set{R}}-\dim\set{B}_0.
\]
\end{lemma}

\section{Main results}\label{sec:main-results}

In this section, we fix a field $\FF\in\{\RR,\CC\}$,  
and let $\tens{T}_1$ and $\tens{T}_2$ be two nonzero
third-order tensors in disjoint factor spaces. 
The first three subsections
concern arbitrary third-order tensors; the remaining subsections treat nonzero symmetric third-order tensors $\tens{S}_1$ and $\tens{S}_2$.  %Since zero summands are trivially additive, they are omitted whenever both summands are assumed to be nonzero.

\subsection{Strassen CP rank gaps controlled by CP rank excess}
 In this subsection, after interchanging the two summands if necessary, we
assume that
\begin{equation}\label{eq:ordinary-excess-order}
e_r(\tens{T}_1)\le e_r(\tens{T}_2).
\end{equation}
For $e\in\mathbb Z_{\ge0}$, define
\begin{equation}\label{eq:cre}
c_{re}(e)\eqdef
\begin{cases}
0,
&0\le e\le4,\\[2pt]
e-4,
&5\le e\le20,\\[2pt]
e-\left\lceil
\dfrac{\sqrt{4e+1}-1}{2}
\right\rceil,
&e\ge21.
\end{cases}
\end{equation}
The function $c_{re}$ is nondecreasing on the nonnegative integers.  The
following theorem uses $c_{re}$ to bound the Strassen CP rank gap in terms of the smaller
CP rank excess.

\begin{theorem}[CP rank excess bound]\label{theorem-1-3}
Under assumption \eqref{eq:ordinary-excess-order}, one has
\begin{equation*}
0\le\Delta_r(\tens{T}_1,\tens{T}_2)
\le c_{re}\bigl(e_r(\tens{T}_1)\bigr).
\end{equation*}
As a direct consequence, we have 
\begin{equation}\label{eq:ordinary-small-excess-additivity}
e_r(\tens{T}_1)\le4
\ \Longrightarrow\ 
\Delta_r(\tens{T}_1,\tens{T}_2)=0.
\end{equation}
\end{theorem}

The proof of \cref{theorem-1-3} relies on the following two constraints
relating a positive gap to the smaller CP rank excess.
The first constraint of \cref{prop:ordinary-residual-excess} gives
\eqref{eq:ordinary-small-excess-additivity}.  The two
constraints together yield the control function $c_{re}$.

\begin{proposition}[CP rank excess constraints]
\label{prop:ordinary-residual-excess}
Suppose that $\Delta_r(\tens{T}_1,\tens{T}_2)>0$.  Then
\begin{equation}\label{eq:ordinary-residual-constraints}
\Delta_r(\tens{T}_1,\tens{T}_2)
\le e_r(\tens{T}_1)-4,
\ \ 
\Delta_r(\tens{T}_1,\tens{T}_2)
\le
\bigl(e_r(\tens{T}_1)-\Delta_r(\tens{T}_1,\tens{T}_2)\bigr)^2.
\end{equation}
\end{proposition}

We next convert the excess estimate into a bound that depends only on
CP rank.  For $k\ge1$, take the supremum over all pairs of third-order
tensors $(\tens{T}_1,\tens{T}_2)$ over $\FF$, with arbitrary formats, and
define the maximal gap at CP rank $k$ as follows:
\begin{equation*}
D_{k,\FF}
\eqdef
\sup\left\{
\Delta_r(\tens{T}_1,\tens{T}_2):
r(\tens{T}_1)=k
\right\},
\ \ 
D_{0,\FF}\eqdef0.
\end{equation*}
For $k\ge1$, set
$z(k)
\eqdef
\left\lceil(\sqrt{8k+1}-1)/2\right\rceil.$ 
The maximum rank estimates of \cite[Theorem~1(1) and Theorem~3]{sumi2010maximal}, combined with
\cref{theorem-1-3}, give the following rank-dependent consequence.

\begin{corollary}[Rank-dependent gap bound]
\label{cor:ordinary-rank-bound}
For every $k\ge1$, one has
\begin{equation}\label{eq:ordinary-rank-capacity-bound}
D_{k,\FF}
\le
c_{re}\bigl(k-z(k)\bigr).
\end{equation}
Therefore, the following estimate holds:
\begin{equation}\label{eq:ordinary-two-rank-bound}
\Delta_r(\tens{T}_1,\tens{T}_2)
\le
\min\left\{
c_{re}\bigl(r(\tens{T}_1)-z(r(\tens{T}_1))\bigr),
c_{re}\bigl(r(\tens{T}_2)-z(r(\tens{T}_2))\bigr)
\right\}.
\end{equation}
This gives the following vanishing range:
\begin{equation}\label{eq:ordinary-rank-eight-capacity}
D_{k,\FF}=0,
\ \  0\le k\le8.
\end{equation}
\end{corollary}

The vanishing range gives a criterion that is symmetric in the two
summands.

\begin{corollary}[Low rank additivity]
\label{cor:ordinary-low-rank-additivity}
We have 
\begin{equation}\label{eq:ordinary-rank-eight-additivity}
\min\{r(\tens{T}_1),r(\tens{T}_2)\}\le8 \ \Longrightarrow\ 
\Delta_r(\tens{T}_1,\tens{T}_2)=0.
\end{equation}
Therefore, the following total rank condition is sufficient:
\begin{equation}\label{eq:ordinary-total-rank-additivity}
r(\tens{T}_1)+r(\tens{T}_2)\le17
\ \Longrightarrow\ 
\Delta_r(\tens{T}_1,\tens{T}_2)=0.
\end{equation}
\end{corollary}

\subsection{Strassen CP rank gaps controlled by the maximum mode rank}
In this subsection, after interchanging the two summands if necessary, we
assume that
\begin{equation}\label{eq:ordinary-dimension-order}
d(\tens{T}_1)\le d(\tens{T}_2).
\end{equation}
The following theorem converts the CP rank excess estimate into a bound
that depends only on the smaller maximum mode rank.

\begin{theorem}[Maximum mode rank gap bound]
\label{thm:ordinary-dimension}
Under assumption \eqref{eq:ordinary-dimension-order}, one has
\begin{equation*}
0\le\Delta_r(\tens{T}_1,\tens{T}_2)
\le
c_{re}\left(
\frac{d(\tens{T}_1)(d(\tens{T}_1)-1)}{2}
\right).
\end{equation*}
Equivalently, the bound can be written as follows:
\begin{equation}\label{eq:ordinary-dimension-piecewise}
0\le\Delta_r(\tens{T}_1,\tens{T}_2)\le
\begin{cases}
0,
&1\le d(\tens{T}_1)\le3,\\[2pt]
\dfrac{d(\tens{T}_1)(d(\tens{T}_1)-1)}{2}-4,
&4\le d(\tens{T}_1)\le6,\\[7pt]
\dfrac{d(\tens{T}_1)(d(\tens{T}_1)-1)}{2}
-\left\lceil
\dfrac{
\sqrt{2d(\tens{T}_1)(d(\tens{T}_1)-1)+1}-1
}{2}
\right\rceil,
&d(\tens{T}_1)\ge7.
\end{cases}
\end{equation}
The first branch gives the following additivity condition:
\begin{equation}\label{eq:ordinary-low-dimension-additivity}
\min\{d(\tens{T}_1),d(\tens{T}_2)\}\le3
\ \Longrightarrow\ 
\Delta_r(\tens{T}_1,\tens{T}_2)=0.
\end{equation}
\end{theorem}

The proof of \cref{thm:ordinary-dimension} requires the following maximum rank estimates.
The first inequality of \cref{prop:ordinary-mode-dimension} follows by applying
\cite[Theorem~1(1) and Theorem~3]{sumi2010maximal} after embedding the mode supports
in spaces of the common dimension $d(\tens{T})$; the second
follows from $e_r(\tens{T})=r(\tens{T})-d(\tens{T})$. 

\begin{proposition}[Maximum rank and excess estimates]
\label{prop:ordinary-mode-dimension}
Every third-order tensor $\tens{T}$ satisfies
\begin{equation*}
r(\tens{T})
\le
\frac{d(\tens{T})(d(\tens{T})+1)}{2},
\ \ 
e_r(\tens{T})
\le
\frac{d(\tens{T})(d(\tens{T})-1)}{2}.
\end{equation*}
\end{proposition}

The above estimate gives a simple quadratic bound with an explicit
constant.

\begin{corollary}[Quadratic maximum mode rank bound]
\label{cor:ordinary-dimension-asymptotics}
Under assumption \eqref{eq:ordinary-dimension-order}, one has
\begin{equation}\label{eq:ordinary-dimension-asymptotic}
\Delta_r(\tens{T}_1,\tens{T}_2)
\le
\frac{d(\tens{T}_1)^2}{2}
-
\left(
\frac12+\frac1{\sqrt2}
\right)d(\tens{T}_1)
+\frac12+\frac1{2\sqrt2}.
\end{equation}
\end{corollary}

\subsection{CP rank additivity criteria and comparison with the literature}
\label{sec:ordinary-additivity-literature}

\Cref{tab:ordinary-additivity} compares the criteria
proved above with the closest conditions in the literature. Unless both summands occur explicitly, a condition is one-sided and the other summand is unrestricted.

\begingroup
\footnotesize
\setlength{\tabcolsep}{4pt}
\renewcommand{\arraystretch}{1.15}
\begin{longtable}{@{}>{\raggedright\arraybackslash}p{.25\linewidth}
>{\centering\arraybackslash}p{.08\linewidth}
>{\raggedright\arraybackslash}p{.19\linewidth}
>{\raggedright\arraybackslash}p{.41\linewidth}@{}}
\caption{Sufficient conditions for CP rank additivity.}
\label{tab:ordinary-additivity}\\
\toprule
Condition & Field & Result in this paper & Related published condition \\
\midrule
\endfirsthead
\caption[]{Sufficient conditions for CP rank additivity (continued).}\\
\toprule
Condition & Field & Result in this paper & Related published condition \\
\midrule
\endhead
\bottomrule
\endfoot

$\min\{e_r(\tens{T}_1),e_r(\tens{T}_2)\}\le4$
& $\RR,\CC$
& \cref{theorem-1-3}, \eqref{eq:ordinary-small-excess-additivity}
& $r(\tens{T})\le d_i(\tens{T})+2$ for one summand $\tens{T}$ and
some mode $i$
\cite[Theorem~1.1]{buczynski2020additivity}.\\[4pt]
\midrule

$\min\{r(\tens{T}_1),r(\tens{T}_2)\}\le8$
& $\RR,\CC$
& \cref{cor:ordinary-low-rank-additivity},
\eqref{eq:ordinary-rank-eight-additivity}
& $\min\{r(\tens{T}_1),r(\tens{T}_2)\}\le6$
\cite[Theorem~1.1]{buczynski2020additivity}.\\[4pt]
\midrule

$r(\tens{T}_1)+r(\tens{T}_2)\le17$
& $\RR,\CC$
& \cref{cor:ordinary-low-rank-additivity},
\eqref{eq:ordinary-total-rank-additivity}
& Over $\RR$, $r(\tens{T}_1)+r(\tens{T}_2)\le13$
is sufficient by
\cite[Theorem~1.1]{buczynski2020additivity}.
Over $\CC$, $r(\tens{T}_1)+r(\tens{T}_2)\le14$
is sufficient by combining the same theorem with
\cite[Theorem~1.0.3(iii)]{rupniewski2024additivity}.\\[4pt]
\midrule

$\min\{d(\tens{T}_1),d(\tens{T}_2)\}\le3$
& $\RR,\CC$
& \cref{thm:ordinary-dimension},
\eqref{eq:ordinary-low-dimension-additivity}
& $d_i(\tens{T}),d_j(\tens{T})\le3$ for one summand $\tens{T}$ and
two distinct modes $i,j$
\cite[Theorem~1.1]{buczynski2020additivity}.\\
\end{longtable}
\endgroup

\subsection{Strassen symmetric rank gaps controlled by symmetric rank excess}
In this subsection,
after interchanging the two summands if necessary, we assume that
\begin{equation}\label{eq:symmetric-excess-order}
e_s(\tens{S}_1)\le e_s(\tens{S}_2).
\end{equation}
The next theorem gives a piecewise bound for the symmetric rank gap and
four sufficient conditions for additivity.

\begin{theorem}[Symmetric rank excess bound]\label{thm:sym-q}
Under assumption \eqref{eq:symmetric-excess-order}, one has
\begin{equation}\label{eq:sym-q-main}
0\le \Delta_s(\tens{S}_1,\tens{S}_2)\le
\begin{cases}
\min\bigl\{
2e_s(\tens{S}_1),
(e_s(\tens{S}_2)-5)_+,
(\Delta_c(\tens{S}_2)-1)_+
\bigr\},
&e_s(\tens{S}_1)\le4,\\[2pt]
\min\bigl\{
2e_s(\tens{S}_1),
e_s(\tens{S}_1)+e_s(\tens{S}_2)-9,
(\Delta_c(\tens{S}_2)+e_s(\tens{S}_1)-5)_+
\bigr\},
&e_s(\tens{S}_1)\ge5.
\end{cases}
\end{equation}
In addition, $\Delta_s(\tens{S}_1,\tens{S}_2)=0$ whenever any of the
following conditions holds:
\begin{enumerate}
\item\label{item:symmetric-excess-zero}
$e_s(\tens{S}_1)=0$;
\item\label{item:symmetric-excess-rectangle}
$e_s(\tens{S}_1)\le4$ and $e_s(\tens{S}_2)\le5$;
\item\label{item:symmetric-equal-rank-partner}
$e_s(\tens{S}_1)\le5$ and $\Delta_c(\tens{S}_2)=0$;
\item\label{item:symmetric-small-comon-partner}
$e_s(\tens{S}_1)\le4$ and $\Delta_c(\tens{S}_2)\le1$.
\end{enumerate}
\end{theorem}

The above quantitative bound in \cref{thm:sym-q} combines the next two estimates.  The first is
one-sided and retains information about the Comon rank gap of the other
summand.

\begin{proposition}[One-sided bounds]\label{prop:one-sided}
For every ordered pair $(\tens{S}_1,\tens{S}_2)$ of nonzero symmetric
third-order tensors in disjoint variable spaces, one has
\begin{equation}\label{eq:one-sided-piecewise}
0\le \Delta_s(\tens{S}_1,\tens{S}_2)\le
\begin{cases}
\min\bigl\{
2e_s(\tens{S}_1),
(\Delta_c(\tens{S}_2)-1)_+
\bigr\},
&e_s(\tens{S}_1)\le4,\\[2pt]
\min\bigl\{
2e_s(\tens{S}_1),
(\Delta_c(\tens{S}_2)+e_s(\tens{S}_1)-5)_+
\bigr\},
&e_s(\tens{S}_1)\ge5.
\end{cases}
\end{equation}
One also has
\begin{equation}\label{eq:minimum-rank-bound}
\Delta_s(\tens{S}_1,\tens{S}_2)
\le \min\{s(\tens{S}_1),s(\tens{S}_2)\}-1.
\end{equation}
If $r(\tens{S}_2)=s(\tens{S}_2)$, then
\begin{equation}\label{eq:eq-partner}
\Delta_s(\tens{S}_1,\tens{S}_2)
\le (e_s(\tens{S}_1)-5)_+.
\end{equation}
\end{proposition}

The following
single-tensor estimate is also used in the proof of \cref{thm:sym-q}.  

\begin{lemma}[Comon rank gap bound]\label{lem:comon-gap-bound}
Every symmetric third-order tensor $\tens{S}$ satisfies
\begin{equation}\label{eq:comon-gap-excess-bound}
0\le\Delta_c(\tens{S})\le\bigl(e_s(\tens{S})-4\bigr)_+.
\end{equation}
Thus the following consequence holds:
\begin{align}
s(\tens{S})&\le9
\ \Longrightarrow\  \Delta_c(\tens{S})=0.
\label{eq:rank-nine-comon-zero}
\end{align}
Moreover, $s(\tens{S})\le10$ implies $\Delta_c(\tens{S})\le1$. 
\end{lemma}

Combining this lemma with \cref{prop:one-sided} yields the bound that
uses both excesses.

\begin{proposition}[Joint excess bound]\label{prop:joint-excess}
Under assumption \eqref{eq:symmetric-excess-order}, one has
\begin{equation}\label{eq:two-sided-excess}
\Delta_s(\tens{S}_1,\tens{S}_2)
\le\max\bigl\{0,
e_s(\tens{S}_1)+e_s(\tens{S}_2)-9,
e_s(\tens{S}_2)-5\bigr\}.
\end{equation}
\end{proposition}

Combining \cref{prop:one-sided} with
\cref{lem:comon-gap-bound} gives the following additivity criterion
for symmetric ranks at most nine and ten.

\begin{corollary}[Additivity for symmetric ranks at most nine and ten]
\label{cor:symmetric-rank-nine-ten}
We have 
\begin{equation*}
\min\{s(\tens{S}_1),s(\tens{S}_2)\}\le9,
\ \ 
\max\{s(\tens{S}_1),s(\tens{S}_2)\}\le10 \ \Longrightarrow\ 
\Delta_s(\tens{S}_1,\tens{S}_2)=0.
\end{equation*}
\end{corollary}

\subsection{Strassen symmetric rank gaps controlled by CP rank excess}
In this subsection,
after interchanging the two summands if necessary, we assume that
\begin{equation}\label{eq:smaller-cp-excess}
e_r(\tens{S}_1)\le e_r(\tens{S}_2).
\end{equation}
The next result connects CP rank excess directly to the
Strassen symmetric rank gap. 

\begin{theorem}[CP rank excess control of the symmetric gap]
\label{thm:sym-p}
Under assumption \eqref{eq:smaller-cp-excess}, one has 
\begin{equation}\label{eq:symmetric-gap-from-cp-excess}
0\le\Delta_s(\tens{S}_1,\tens{S}_2)\le
\begin{cases}
2e_r(\tens{S}_1),&e_r(\tens{S}_1)\le3,\\[2pt]
12e_r(\tens{S}_1)-9,&e_r(\tens{S}_1)\ge4.
\end{cases}
\end{equation}
Moreover, the following condition implies additivity:
\begin{equation}
e_r(\tens{S}_2)\le3
\ \Longrightarrow\ 
\Delta_s(\tens{S}_1,\tens{S}_2)=0.
\label{eq:cp-excess-three-symmetric-additivity}
\end{equation}
\end{theorem}

\subsection{Strassen symmetric rank gaps controlled by essential dimension}

In this subsection, after interchanging the two
summands if necessary, we assume  that
\begin{equation}\label{eq:dimension-order}
d(\tens{S}_1)\le d(\tens{S}_2).
\end{equation}
The next
result bounds the gap using only the smaller essential dimension.

\begin{theorem}[Dimension-dependent gap bound]\label{thm:sym-d}
Under assumption \eqref{eq:dimension-order}, one has
\begin{equation}\label{eq:sym-d-main}
0\le \Delta_s(\tens{S}_1,\tens{S}_2)\le
\begin{cases}
2d(\tens{S}_1)-2,
&1\le d(\tens{S}_1)\le4,\\[2pt]
16,
&d(\tens{S}_1)=6,\\[2pt]
2\left\lceil
\dfrac{d(\tens{S}_1)(d(\tens{S}_1)+1)}{6}
\right\rceil,
&d(\tens{S}_1)\ge5,\  d(\tens{S}_1)\ne6.
\end{cases}
\end{equation}
Moreover, one has
\begin{align}
\min\{d(\tens{S}_1),d(\tens{S}_2)\}\le1
&\ \Longrightarrow\ 
\Delta_s(\tens{S}_1,\tens{S}_2)=0,
\label{eq:symmetric-dimension-one-additivity}\\
\max\{d(\tens{S}_1),d(\tens{S}_2)\}\le4
&\ \Longrightarrow\ 
\Delta_s(\tens{S}_1,\tens{S}_2)=0.
\label{eq:symmetric-dimension-four-additivity}
\end{align}
\end{theorem}

The piecewise estimate also gives a simpler bound.

\begin{corollary}[Uniform dimension bound]
\label{cor:symmetric-dimension-asymptotics}
Under assumption \eqref{eq:dimension-order}, one has
\begin{equation}\label{eq:symmetric-uniform-dimension}
\Delta_s(\tens{S}_1,\tens{S}_2)
\le
\frac{d(\tens{S}_1)^2}{3}
+\frac{d(\tens{S}_1)}{3}
+2.
\end{equation}
\end{corollary}

The proof of \cref{thm:sym-d} uses the formula for $g_3$ in
\eqref{eq:generic} and the following shifted maximum rank estimate. 
The real case of the shifted estimate is
\cite[Proposition~6.6]{li2026comonrankgap}.
A proof is given in \cref{sec:shifted-proof}.

\begin{proposition}[Shifted maximum symmetric rank bound]\label{prop:shifted}
For $n\ge2$, one has
\[
M_{3,\FF}(n)\le2g_3(n-1)+1.
\]
For $1\le n\le4$, the estimate
$M_{3,\FF}(n)\le2n-1$ holds.  Consequently, under assumption
\eqref{eq:dimension-order}, the gap satisfies
\begin{equation*}
\Delta_s(\tens{S}_1,\tens{S}_2)
\le M_{3,\FF}(d(\tens{S}_1))-1.
\end{equation*}
\end{proposition}

\subsection{Sufficient additivity criteria and comparison with the literature}
\label{sec:quantitative-literature}

\Cref{tab:symmetric-additivity} compares the main
criteria proved here with related results in the literature.
No ordering of the summands is assumed. Conditions involving
$\tens{S}_i,\tens{S}_j$ require that there exist indices with $\{i,j\}=\{1,2\}$ satisfying the stated conditions.

\begingroup
\makeatletter
\long\def\LT@makecaption#1#2#3{%
  \LT@mcol\LT@cols c{\hbox to\z@{\hss
    \parbox[t]\linewidth{%
      \@tempcnta=64\relax
      \@makecaption{#1{#2}}{#3}}%
    \hss}}}
\makeatother
\footnotesize
\setlength{\tabcolsep}{4pt}
\renewcommand{\arraystretch}{1.15}
\begin{longtable}{@{}>{\raggedright\arraybackslash}p{.23\linewidth}
>{\centering\arraybackslash}p{.06\linewidth}
>{\raggedright\arraybackslash}p{.20\linewidth}
>{\raggedright\arraybackslash}p{.44\linewidth}@{}}
\caption{Sufficient conditions for symmetric rank additivity.}
\label{tab:symmetric-additivity}\\
\toprule
Condition & Field & Result in this paper & Related published condition \\
\midrule
\endfirsthead
\caption[]{Sufficient conditions for symmetric rank additivity
(continued).}\\
\toprule
Condition & Field & Result in this paper & Related published condition \\
\midrule
\endhead
\endfoot
\bottomrule
\endlastfoot

$\min\{e_s(\tens{S}_1),e_s(\tens{S}_2)\}=0$
& $\RR,\CC$
& \cref{thm:sym-q} \ref{item:symmetric-excess-zero}
& The same condition over $\CC$
\cite[Corollary~3.2]{carlini2015symmetricstrassen}.\\[4pt]
\midrule

$\min\{e_s(\tens{S}_1),e_s(\tens{S}_2)\}\le4$,
$\max\{e_s(\tens{S}_1),e_s(\tens{S}_2)\}\le5$
& $\RR,\CC$
& \cref{thm:sym-q} \ref{item:symmetric-excess-rectangle}
& $\max\{e_s(\tens{S}_1),e_s(\tens{S}_2)\}\le2$
\cite[Theorem~1.1]{friedland2016remarks}
\cite[Theorem~1.1]{buczynski2020additivity}.\\[4pt]
\midrule

$e_s(\tens{S}_i)\le5$, $\Delta_c(\tens{S}_j)=0$
& $\RR,\CC$
& \cref{prop:one-sided}, \eqref{eq:eq-partner}
& $e_s(\tens{S}_i)\le2$ and $r(\tens{S}_j)=s(\tens{S}_j)$
\cite[Theorem~1.1]{friedland2016remarks}
\cite[Theorem~1.1]{buczynski2020additivity}.\\[4pt]
\midrule

$e_s(\tens{S}_i)\le4$, $\Delta_c(\tens{S}_j)\le1$
& $\RR,\CC$
& \cref{prop:one-sided}, \eqref{eq:one-sided-piecewise}
& \qquad\qquad\qquad\qquad ------
\\[4pt]
\midrule

$\max\{e_r(\tens{S}_1),e_r(\tens{S}_2)\}\le3$
& $\RR,\CC$
& \cref{thm:sym-p}, \eqref{eq:cp-excess-three-symmetric-additivity}
& $\max\{e_r(\tens{S}_1),e_r(\tens{S}_2)\}\le1$,
obtained by combining
\cite[Theorem~1.1]{friedland2016remarks} with
\cite[Theorem~1.1]{buczynski2020additivity}.
\\[4pt]
\midrule

$\min\{s(\tens{S}_1),s(\tens{S}_2)\}\le9$,\newline
$\max\{s(\tens{S}_1),s(\tens{S}_2)\}\le10$
& $\RR,\CC$
& \cref{cor:symmetric-rank-nine-ten}
& Over $\CC$, $\max\{s(\tens{S}_1),s(\tens{S}_2)\}\le7$
\cite[Corollary~1.4]{seigal2020ranks}
\cite[Theorem~1.0.3(iii)]{rupniewski2024additivity}.\\[4pt]
\midrule

$\min\{d(\tens{S}_1),d(\tens{S}_2)\}\le1$
& $\RR,\CC$
& \cref{thm:sym-d}, \eqref{eq:symmetric-dimension-one-additivity}
& A pure power summand over $\CC$
\cite[Proposition~3.1]{carlini2015symmetricstrassen}.\\[4pt]
\midrule

$\max\{d(\tens{S}_1),d(\tens{S}_2)\}\le4$
& $\RR,\CC$
& \cref{thm:sym-d}, \eqref{eq:symmetric-dimension-four-additivity}
& Over $\CC$, the same condition follows by combining the bound
$r=s\le7$ for each summand from
\cite[Section~2, Theorem~1.3]{seigal2020ranks}
with \cite[Theorem~1.0.3]{rupniewski2024additivity}.
Over $\RR$, $\max\{d(\tens{S}_1),d(\tens{S}_2)\}\le3$ from
\cite[Proposition~4.2]{seigal2020ranks} and
\cite[Theorem~1.1]{buczynski2020additivity}.\\
\end{longtable}
\endgroup

\section{Proofs of the CP rank bounds}\label{sec:ordinary-proofs}

We first prove the CP rank excess constraints in
\cref{prop:ordinary-residual-excess}. 
A matrix space $\set{W}\subseteq\set{B}\otimes\set{C}$
has hook shape $(u,v)$ if there exist subspaces
$\set{E}\subseteq\set{B}$ and $\set{F}\subseteq\set{C}$
such that
\[
\set{W}\subseteq\set{E}\otimes\set{C}+\set{B}\otimes\set{F},
\ \ 
\dim\set{E}\le u,
\ 
\dim\set{F}\le v.
\]
Two additivity criteria in the literature will be used below.  A tensor is
CP rank additive with every partner if one of its slice spaces has
hook shape $(1,2)$, up to transposition
\cite[Theorem~1.0.3(i)]{rupniewski2024additivity}.  The same conclusion
holds when two of its actual mode dimensions are at most three
\cite[Theorem~1.1]{buczynski2020additivity}.

\subsection{Projection kernels and the quadratic constraint}

Assume that $\Delta_r(\tens{T}_1,\tens{T}_2)>0$.  Choose a mode of
$\tens{T}_1$ whose actual dimension is $d(\tens{T}_1)$, and take the slice
spaces of both tensors in this mode.  Write
$\set{W}_1\subseteq\set{B}_1\otimes\set{C}_1$ and 
$\set{W}_2\subseteq\set{B}_2\otimes\set{C}_2$, 
where the four matrix factors are the actual supports, so that
$\dim\set{W}_1=d(\tens{T}_1)$.  
Let $\set{R}\subseteq
(\set{B}_1\oplus\set{B}_2)\otimes
(\set{C}_1\oplus\set{C}_2)$
have minimum dimension among spaces spanned by rank-one matrices
and containing $\set{W}_1\oplus\set{W}_2$.
Then
\begin{equation}\label{eq:ordinary-proof-residual}
\dim\set{R}-\dim\set{W}_1-r(\set{W}_2)
=e_r(\tens{T}_1)-\Delta_r(\tens{T}_1,\tens{T}_2).
\end{equation}
We use the expression in slice spaces on the left throughout the
reductions below, including when the chosen mode is no longer maximal.

Apply the projection onto $\set{C}_2$ along $\set{C}_1$
in the second tensor factor of $\set{R}$.
The resulting space has first-factor support
$\set{B}_2\oplus\set{E}_1$.
Similarly, projection onto $\set{B}_2$ along $\set{B}_1$
in the first tensor factor gives a space with second-factor support
$\set{C}_2\oplus\set{F}_1$.
Thus $\set{E}_1\subseteq\set{B}_1$ and
$\set{F}_1\subseteq\set{C}_1$ describe the extra support
of these projections.
Interchanging the two summands defines
$\set{E}_2\subseteq\set{B}_2$ and
$\set{F}_2\subseteq\set{C}_2$.
Put 
\[
\lambda_1\eqdef\dim\set{E}_1,
\  \mu_1\eqdef\dim\set{F}_1,
\ \ 
\lambda_2\eqdef\dim\set{E}_2,
\  \mu_2\eqdef\dim\set{F}_2.
\]
These spaces are the leakage spaces of the chosen covering space $\set{R}$.  The
following decomposition records the possible supports of its
basis matrices of rank one.

\begin{proposition}[An adapted decomposition into seven classes]
\label{lem:ordinary-seven-class}
A basis of $\set{R}$ consisting of rank-one matrices can be partitioned as
\[
\set{P}_1\sqcup\set{P}_2\sqcup
\set{H}_1\sqcup\set{H}_2\sqcup
\set{V}_1\sqcup\set{V}_2\sqcup\set{M}.
\]
Among all such bases and partitions satisfying the containment conditions below, first maximize
$|\set{P}_1|+|\set{P}_2|$ and then, subject to this maximum, maximize
$|\set{H}_1|+|\set{H}_2|+|\set{V}_1|+|\set{V}_2|$.
Each class is contained in the tensor product space listed below: 
\[
\begin{array}{c|c}
\text{class} & \text{containing tensor product space}\\
\hline
\set{P}_1 & \set{B}_1\otimes\set{C}_1\\
\set{P}_2   & \set{B}_2\otimes\set{C}_2\\
\set{H}_1    & \set{E}_1\otimes(\set{C}_1\oplus\set{F}_2)\\
\set{H}_2    & \set{E}_2\otimes(\set{F}_1\oplus\set{C}_2)\\
\set{V}_1    & (\set{B}_1\oplus\set{E}_2)\otimes\set{F}_1\\
\set{V}_2    & (\set{E}_1\oplus\set{B}_2)\otimes\set{F}_2\\
\set{M}   &
(\set{E}_1\oplus\set{E}_2)\otimes
(\set{F}_1\oplus\set{F}_2).
\end{array}
\]
Let $p_1,p_2,h_1,h_2,v_1,v_2,m$ be their respective cardinalities.
Then
\[
\dim\set{R}=p_1+p_2+h_1+h_2+v_1+v_2+m
\]
and
\begin{align}
r(\set{W}_1)
&\le p_1+h_1+v_1+\min\{m,\lambda_1\mu_1\},
\label{eq:ordinary-seven-class-T}\\
r(\set{W}_2)
&\le p_2+h_2+v_2+\min\{m,\lambda_2\mu_2\}.
\label{eq:ordinary-seven-class-S}
\end{align}
\end{proposition}

For the support classification and basis convention, see
\cite[Lemma~6.1.1 and Notation~6.1.2]{rupniewski2024additivity}.
The rank inequalities follow from
\cite[Proposition~6.1.3]{rupniewski2024additivity}. 
The next proposition bounds the leakage dimensions and the gap
using the residual for the slice spaces.

\begin{proposition}[Leakage bounds]\label{lem:ordinary-leakage}
Let $\set{W}_1$ and $\set{W}_2$ be slice spaces in their actual matrix
supports, and let $\set{R}$ have minimum dimension among spaces spanned by
matrices of rank one and containing their direct sum.  Assume that
$r(\set{W}_1)+r(\set{W}_2)-\dim\set{R}>0.$ 
Then the leakage dimensions satisfy
\begin{equation*}
\begin{aligned}
&\lambda_1\le\dim\set{R}-\dim\set{W}_1-r(\set{W}_2),\ \
\mu_1\le\dim\set{R}-\dim\set{W}_1-r(\set{W}_2),\\
&r(\set{W}_1)+r(\set{W}_2)-\dim\set{R}
\le\lambda_1\mu_1
\le\bigl(\dim\set{R}-\dim\set{W}_1-r(\set{W}_2)\bigr)^2.
\end{aligned}
\end{equation*}
In particular, both $\lambda_1$ and $\mu_1$ are positive.
\end{proposition}

\begin{proof}
The first two inequalities are
\cite[Lemma~3.5]{buczynski2020additivity}.
To obtain the product bound, add
\eqref{eq:ordinary-seven-class-T} and
\eqref{eq:ordinary-seven-class-S}.  Since the second minimum is at
most $m$, this gives
\[
r(\set{W}_1)+r(\set{W}_2)
\le\dim\set{R}+\min\{m,\lambda_1\mu_1\}.
\]
This proves the product bound.  The first two inequalities give
the squared bound, and the positive gap makes both leakage dimensions
positive.
\end{proof}

% \begin{remark}[Relation to projection bounds]
% The projection estimates in
% \cite[Lemma~3.5]{buczynski2020additivity} bound the two leakage
% dimensions separately.  Combining them with the rank inequalities for the seven classes gives the quadratic constraint in
% \eqref{eq:ordinary-leakage}.  Thus this part of the argument extracts
% a quantitative gap bound from the published estimates for matrix spaces.
% \end{remark}

\subsection{A common hook and the linear constraint}

We first remove the pure terms on the $\tens{T}_1$-side by repletion and
digestion.  Suppose that a minimum cover of length $N$ contains
$\matr{X}\in\set{P}_1$, and define
$\varepsilon\eqdef
\dim(\set{W}_1+\Span\{\matr{X}\})-\dim\set{W}_1\in\{0,1\}.$ 
By \cite[Corollary~6.3.5]{rupniewski2024additivity}, the
reduction leaves $\set{W}_2$ unchanged, does not decrease the gap,
and changes the cover length and slice dimension as follows:
\[
(N,\dim\set{W}_1)
\longmapsto
(N-1,\dim\set{W}_1-1+\varepsilon).
\]
Thus the residual $\dim\set{R}-\dim\set{W}_1-r(\set{W}_2)$ decreases
by $\varepsilon$.
After each step, restrict to the actual supports using
\cite[Proposition~2.10(i)]{buczynski2020additivity} and choose an
adapted minimum cover.  Since $N$ decreases, the process terminates
with no $\set{P}_1$ terms.  The resulting pair still has positive
gap and its residual is no larger than the original one. 
For this reduced pair, a single projection kernel controls a common
hook for the entire slice space.

\begin{proposition}[Common hook reduction]\label{lem:ordinary-common-hook}
Assume that the matrix factors are the actual supports and that the
adapted minimum basis satisfies $\set{P}_1=\varnothing$.
If $\set{E}_1\subsetneq\set{B}_1$, then there is a subspace
$\set{F}_*\subseteq\set{C}_1$ such that
\begin{equation}\label{eq:ordinary-common-hook}
\begin{gathered}
\set{W}_1\subseteq
\set{E}_1\otimes\set{C}_1+\set{B}_1\otimes\set{F}_*,\ \
\dim\set{F}_*\le
\dim\set{R}-\dim\set{W}_1-r(\set{W}_2)-\lambda_1.
\end{gathered}
\end{equation}
The transposed statement holds when
$\set{F}_1\subsetneq\set{C}_1$.
\end{proposition}

\begin{proof}
Let $\rho_{\set{C}_2}$ denote the projection in the second tensor factor onto $\set{C}_2$ along $\set{C}_1$, and set
$\set{K}_C\eqdef\ker(\rho_{\set{C}_2}|_{\set{R}})$.  
The support in the first factor of $\rho_{\set{C}_2}(\set{R})$ is
$\set{B}_2\oplus\set{E}_1$.  Hence
\cref{lem:ambient-support} gives
$\dim\rho_{\set{C}_2}(\set{R})\ge r(\set{W}_2)+\lambda_1.$ 
Since $\set{W}_1\subseteq\set{K}_C$, we have
$$\dim(\set{K}_C/\set{W}_1)
\le\dim\set{R}-\dim\set{W}_1-r(\set{W}_2)-\lambda_1.$$
Write each $\set{V}_1$ basis element as
\[
\matr{X}_i
=(\vect{b}_i^{(1)}+\vect{b}_i^{(2)})\otimes\vect{c}_i^{(1)},
\ \ 
\vect{c}_i^{(1)}\in\set{F}_1.
\]
Each $\matr{X}_i$ belongs to $\set{K}_C$.
Moreover, $\vect{b}_i^{(2)}\ne0$.
Otherwise, $\matr{X}_i$ would lie in $\set{B}_1\otimes\set{C}_1$ and, by the maximality of the pure classes, would belong to $\set{P}_1$, contradicting
$\set{P}_1=\varnothing$. 
Define the space of their lower left blocks by
$\set{L}_C\eqdef
\Span\{\vect{b}_i^{(2)}\otimes\vect{c}_i^{(1)}:
\matr{X}_i\in\set{V}_1\}.$ 
Projection onto the lower left block vanishes on $\set{W}_1$, so
$\dim\set{L}_C\le\dim(\set{K}_C/\set{W}_1)$.
Choose a basis of $\set{L}_C$ from its displayed generators and let
$\set{F}_*$ be the span of their second factors.  Then
\[
\dim\set{F}_*\le\dim\set{L}_C
\le\dim(\set{K}_C/\set{W}_1)
\le\dim\set{R}-\dim\set{W}_1-r(\set{W}_2)-\lambda_1.
\]
Every remaining lower left block belongs to
$\set{B}_2\otimes\set{F}_*$.  Since its first factor is nonzero,
its second factor also belongs to $\set{F}_*$.  Thus every
$\set{V}_1$ term has its second factor in $\set{F}_*$.

Finally, project the $\tens{T}_1$-diagonal block onto
$(\set{B}_1/\set{E}_1)\otimes\set{C}_1$.
The three classes $\set{P}_2,\set{H}_2,\set{V}_2$ have
zero $\tens{T}_1$-diagonal block, while $\set{H}_1$ and $\set{M}$
vanish modulo $\set{E}_1$.  There are no $\set{P}_1$ terms.
Only the $\set{V}_1$ terms remain, so their common space of second factors $\set{F}_*$ gives \eqref{eq:ordinary-common-hook}.
Interchanging the two matrix factors proves the transposed statement.
\end{proof}

% \begin{remark}[Relation to additivity for hook shapes]
% Repletion and digestion are the reductions in
% \cite[Corollary~6.3.5]{rupniewski2024additivity}; additivity for
% spaces of hook shape $(1,2)$ is
% \cite[Theorem~1.0.3(i)]{rupniewski2024additivity}.
% The additional estimate in \cref{lem:ordinary-common-hook} bounds the
% sum of the dimensions of the two hook arms by
% $\dim\set{R}-\dim\set{W}_1-r(\set{W}_2)$.
% It allows the published additivity criterion to exclude positive gaps
% with residual at most three.
% \end{remark}

We can now combine the projection and common hook estimates.

\begin{proof}[Proof of \cref{prop:ordinary-residual-excess}]
Assume that $\Delta_r(\tens{T}_1,\tens{T}_2)>0$.
Applied to the original slice pair,
\cref{lem:ordinary-leakage} and \eqref{eq:ordinary-proof-residual} give
\[
\Delta_r(\tens{T}_1,\tens{T}_2)
\le\bigl(e_r(\tens{T}_1)-\Delta_r(\tens{T}_1,\tens{T}_2)\bigr)^2.
\]
Suppose that $e_r(\tens{T}_1)-\Delta_r(\tens{T}_1,\tens{T}_2)\le3$.
Repletion and digestion give a pair in its actual supports with no
$\set{P}_1$ terms, positive gap, and
\[
\dim\set{R}-\dim\set{W}_1-r(\set{W}_2)
\le e_r(\tens{T}_1)-\Delta_r(\tens{T}_1,\tens{T}_2)\le3,
\]
where the slice spaces and cover on the left are those of the reduced
pair, while $\tens{T}_1$ and $\tens{T}_2$ on the right are the original
tensors.  By \cref{lem:ordinary-leakage}, both $\lambda_1$ and
$\mu_1$ are positive. 
If $\set{E}_1\subsetneq\set{B}_1$, then
\cref{lem:ordinary-common-hook} places $\set{W}_1$ in a hook whose
arm dimensions have sum at most three.  
The second arm is nonzero: otherwise,
$\set{W}_1\subseteq\set{E}_1\otimes\set{C}_1$,
contradicting the fact that $\set{B}_1$ is the actual first-factor
support and $\set{E}_1\subsetneq\set{B}_1$.
The first arm has positive dimension $\lambda_1$.
Since both arm dimensions are positive integers with sum at most
three, $\set{W}_1$ has hook shape $(1,2)$ up to transposition.
By \cite[Theorem~1.0.3(i)]{rupniewski2024additivity},
rank additivity holds, a contradiction.
Therefore
$\set{E}_1=\set{B}_1$.  The transposed argument gives
$\set{F}_1=\set{C}_1$.  Consequently
$\dim\set{B}_1=\lambda_1\le3$ and 
$\dim\set{C}_1=\mu_1\le3.$ 
The criterion for small formats in
\cite[Theorem~1.1]{buczynski2020additivity} again gives
additivity, a contradiction.  Hence
$$e_r(\tens{T}_1)-\Delta_r(\tens{T}_1,\tens{T}_2)\ge4.$$ 
Together with the quadratic constraint proved above, this gives
\eqref{eq:ordinary-residual-constraints}.
\end{proof}

\subsection{CP rank excess, CP rank, and dimension bounds}

Solving the two residual constraints gives the control function
$c_{re}$.  The remaining proofs then use its monotonicity.

\begin{proof}[Proof of \cref{theorem-1-3}]
The assertion is immediate when $\Delta_r(\tens{T}_1,\tens{T}_2)=0$.
Otherwise, put $e\eqdef e_r(\tens{T}_1)$ and
$u\eqdef e-\Delta_r(\tens{T}_1,\tens{T}_2)$.
The two constraints in \cref{prop:ordinary-residual-excess} say that
$u\ge4$ and $e-u\le u^2$.  Solving the latter inequality and using
integrality gives
\[
u\ge\max\left\{4,
\left\lceil\frac{\sqrt{4e+1}-1}{2}\right\rceil\right\}.
\]
For $e\le20$, the ceiling is at most four; for $e\ge21$, it is
at least five.  Thus \eqref{eq:cre} gives
$\Delta_r(\tens{T}_1,\tens{T}_2)\le c_{re}(e)$.
A positive gap also requires $e\ge5$, which proves
\eqref{eq:ordinary-small-excess-additivity}.
\end{proof}

The argument applies to either chosen summand.  Equivalently, symmetry
of the gap and monotonicity of $c_{re}$ give
\begin{equation}\label{eq:ordinary-either-summand}
\Delta_r(\tens{T}_1,\tens{T}_2)
\le\min\{c_{re}(e_r(\tens{T}_1)),c_{re}(e_r(\tens{T}_2))\}.
\end{equation}
We use this form to obtain the following consequences.

\begin{proof}[Proof of \cref{cor:ordinary-rank-bound}]
Consider any tensor $\tens{T}$ of CP rank $k\ge1$.
The maximum CP rank estimate in \cref{prop:ordinary-mode-dimension} gives
$k\le d(\tens{T})(d(\tens{T})+1)/2.$
Thus $d(\tens{T})\ge z(k)$ and
$e_r(\tens{T})\le k-z(k)$.  
By \eqref{eq:ordinary-either-summand} and monotonicity of $c_{re}$, the gap is at most $c_{re}(k-z(k))$ for every tensor of CP rank $k$ and every partner.
Taking the supremum over all such tensors and partners proves \eqref{eq:ordinary-rank-capacity-bound}. 
Applying the same estimate to each summand gives
\eqref{eq:ordinary-two-rank-bound}.
For $1\le k\le8$, the inequality $k-z(k)\le4$ gives
$D_{k,\FF}=0$; the case $k=0$ holds by definition.
\end{proof}

\begin{proof}[Proof of \cref{cor:ordinary-low-rank-additivity}]
The first assertion is \eqref{eq:ordinary-rank-eight-capacity}
applied to the summand of smaller CP rank.
If both CP ranks were at least nine, their sum would be at least
eighteen.  Therefore the total CP rank condition implies the first
assertion as well.
\end{proof}

\begin{proof}[Proof of \cref{thm:ordinary-dimension}]
By \eqref{eq:ordinary-either-summand},
\cref{prop:ordinary-mode-dimension}, and monotonicity, we obtain
\[
\Delta_r(\tens{T}_1,\tens{T}_2)
\le c_{re}\bigl(e_r(\tens{T}_1)\bigr)
\le c_{re}\left(
\frac{d(\tens{T}_1)(d(\tens{T}_1)-1)}2
\right).
\]
Substituting \eqref{eq:cre} gives
\eqref{eq:ordinary-dimension-piecewise}, whose first branch proves
\eqref{eq:ordinary-low-dimension-additivity}.
\end{proof}

\begin{proof}[Proof of \cref{cor:ordinary-dimension-asymptotics}]
For $d(\tens{T}_1)\ge7$, use $\lceil x\rceil\ge x$ and the inequality
\[
\sqrt{2d(\tens{T}_1)(d(\tens{T}_1)-1)+1}
\ge\sqrt2\,d(\tens{T}_1)-\frac1{\sqrt2}
\]
in \eqref{eq:ordinary-dimension-piecewise}.  This gives
\eqref{eq:ordinary-dimension-asymptotic}.  For
$1\le d(\tens{T}_1)\le6$, direct substitution in the first two branches
gives the same bound.
\end{proof}

\section{Polar extraction and regular decompositions}\label{sec:polar-tools}
Throughout this section, let $\set{V}$ be a finite-dimensional
vector space over the fixed field $\FF$.

\subsection{Regular neighborhoods}
For $m\ge1$, define
\[
\Phi_m:\set{V}^m\longrightarrow\Sym^3(\set{V}),
\ \ 
\Phi_m(\vect{\ell}_1,\ldots,\vect{\ell}_m)
\eqdef\sum_{i=1}^m\vect{\ell}_i^3.
\]
We call a tuple \emph{regular} if the differential of $\Phi_m$ is
surjective at that tuple. In terms of normalized symmetric products,
this condition is
$\Span\{\vect{\ell}_i^2\vect{v}:
1\le i\le m,\ \vect{v}\in\set{V}\}=\Sym^3(\set{V}).$ 
The following local property is used in the maximum symmetric rank bound.

\begin{proposition}[Regular cube neighborhoods]
\label{lem:regular-degeneration}
Let $K\eqdef\Phi_m(\vect{\ell}_1,\ldots,\vect{\ell}_m)$ arise from a
regular tuple.  Then $K$ has a Euclidean neighborhood in which every
form admits a regular decomposition into $m$ cubes.
\end{proposition}
\begin{proof}
The set of tuples at which the differential of $\Phi_m$ is
surjective is open.  
Restrict $\Phi_m$ to a
neighborhood of the given tuple on which its differential remains
surjective.  Over $\CC$, regard $\Phi_m$ as a smooth map of real vector
spaces.  The rank theorem \cite[Theorem~4.12]{lee2012smooth}
implies that its image contains a neighborhood of $K$.
Every form in that neighborhood therefore has a preimage that is a
regular tuple.
\end{proof}

% \begin{remark}[Relation to rank arguments using open sets]
% Blekherman and Teitler use open sets of fixed rank to bound maximum
% rank by twice the generic rank
% \cite[Theorems~1--3 and their proofs]{blekherman2015maximum}.
% Here the surjective differential supplies an open set around a
% specified cube sum. This allows us to use the same local principle
% after prescribing a polar and subtracting a controlled number of cubes.
% \end{remark}

\subsection{Polar extraction}
For a nonzero covector $\alpha\in\set{V}^*$, the kernel of
contraction by $\alpha$ on $\Sym^3(\set{V})$ is
$\Sym^3(\ker\alpha)$.
The following construction adapts the polar deformation in
\cite[proof of Proposition~6.6]{li2026comonrankgap} to both fields
and allows the coefficients of the leading cubic to be prescribed.

\begin{proposition}[Polar extraction with a prescribed leading cubic]
\label{lem:polar-leading}
Let $F\in\Sym^3(\set{V})$ and let $\alpha\in\set{V}^*$ satisfy
$F(\alpha)\ne0$. Set
$$Q\eqdef C_F(\alpha),\  
\set{W}\eqdef\Ess(Q),\  
\set{K}\eqdef\set{W}\cap\ker\alpha,\  h\eqdef\rank Q.$$  
Then $\dim\set{K}=h-1$, and $\set{K}$ has a basis
$\vect{w}_1,\ldots,\vect{w}_{h-1}$ with the following property.
For any $c_1,\ldots,c_{h-1}\in\FF\setminus\{0\}$,
there exists a family $S_\varepsilon$ of sums of $h$ cubes
in $\Sym^3(\set{W})$, defined for all sufficiently small
nonzero $\varepsilon\in\FF$, such that 
\[
C_{S_\varepsilon}(\alpha)=Q,
\ \ 
\varepsilon S_\varepsilon\longrightarrow
P\eqdef\sum_{i=1}^{h-1}c_i\vect{w}_i^3
\ \ (\varepsilon\to0).
\]
Consequently, the remainder satisfies
$F-S_\varepsilon\in\Sym^3(\ker\alpha)$.
\end{proposition}
\begin{proof}
Let $b$ be the inverse bilinear form of $Q$ on $\set{W}$, and let
$Q^\sharp:\set{W}^*\to\set{W}$ be the linear map defined by $Q$.
Set
\[
\vect{v}_0\eqdef Q^\sharp(\alpha|_{\set{W}}),
\ \ 
\rho\eqdef\alpha(\vect{v}_0)
=b(\vect{v}_0,\vect{v}_0)=F(\alpha)\ne0.
\]
Therefore, the space $\set{K}=\vect{v}_0^{\perp_b}$ is nondegenerate and has
dimension $h-1$. Choose a $b$-orthogonal basis
$\vect{w}_1,\ldots,\vect{w}_{h-1},\vect{w}_h\eqdef\vect{v}_0$ of
$\set{W}$, and set
$\delta_i\eqdef b(\vect{w}_i,\vect{w}_i)^{-1}.$ 
In this basis, the quadratic form is
$Q=\sum_i\delta_i\vect{w}_i^2$.
For $i<h$, put $k_i\eqdef\delta_i/(\rho c_i)$ and define a linear map
$A$ by
\[
A\vect{w}_i\eqdef k_i\vect{v}_0,
\ \ 
A\vect{v}_0\eqdef
-\sum_{i<h}k_i\rho\delta_i\vect{w}_i.
\]
The map $A$ is skew-adjoint with respect to $b$.
For all sufficiently small $\varepsilon\in\FF$,
$\id-\varepsilon A/2$ is invertible, and the Cayley transform
$G_\varepsilon\eqdef
(\id+\varepsilon A/2)(\id-\varepsilon A/2)^{-1}$
preserves $b$ and its inverse tensor $Q$.
It satisfies
$G_\varepsilon=\id+\varepsilon A+O(\varepsilon^2)$.
For $\vect{v}_i\eqdef G_\varepsilon\vect{w}_i$, it follows that 
\[
Q=\sum_i\delta_i\vect{v}_i^2,
\ \ 
\alpha(\vect{v}_i)=\frac{\varepsilon\delta_i}{c_i}
+O(\varepsilon^2)\ (i<h),
\ \ 
\alpha(\vect{v}_h)=\rho+O(\varepsilon).
\]
For all sufficiently small nonzero $\varepsilon\in\FF$,
the values $\alpha(\vect{v}_i)$ are nonzero for every
$1\le i\le h$.
We may therefore define
$S_\varepsilon\eqdef
\sum_{i=1}^h\frac{\delta_i}{\alpha(\vect{v}_i)}\vect{v}_i^3$. 
Contraction gives $C_{S_\varepsilon}(\alpha)=Q$.
Multiplying by $\varepsilon$ and taking the limit term by term
gives the prescribed $P$. Finally, the contraction kernel identity gives
$F-S_\varepsilon\in\Sym^3(\ker\alpha)$.
\end{proof}

\section{Rank bounds in small dimensions and rank equality}\label{sec:rank-inputs}

\subsection{A bound of seven cubes for real quaternary cubics}
\label{sec:quaternary-rank-bound}
The following bound is \cite[Proposition~5.2]{li2026comonrankgap}.

\begin{proposition}[Real cubics in four variables]
\label{prop:four-variable-seven}
Every real cubic in at most four variables is a sum of at most seven
real cubes. Equivalently, one has $M_{3,\RR}(4)\le7$.
\end{proposition}

The complex bound follows from
\cite[Section~2 and Theorem~1.3]{seigal2020ranks}.
For use below, we record the polarization identities
\begin{align*}
x^2M&=\frac{(x+M)^3-(x-M)^3-2M^3}{6},\\
24UVW&=(U+V+W)^3+(U-V-W)^3+(-U+V-W)^3+(-U-V+W)^3.
\end{align*}

\subsection{Rank equality through CP rank excess three}
\label{sec:cp-excess-three}
The real case of the following rank equality is
\cite[Proposition~3.3]{li2026comonrankgap}.
The wedge support estimates are also needed for the
symmetric rank gap bounds in \cref{sec:gap-proofs}.

\begin{proposition}[Rank equality at CP rank excess at most three]
\label{prop:e3}
Let $\tens{S}$ be a symmetric third-order tensor over either field.
Then
\[
e_r(\tens{S})\le3
\ \Longrightarrow\ 
r(\tens{S})=s(\tens{S}).
\]
The condition $e_s(\tens{S})\le4$ also implies this rank equality.
\end{proposition}

\subsubsection{Wedge supports}\label{sec:wedge-control}
A \emph{minimum partially symmetric CP decomposition} is a minimum CP decomposition
that, after a permutation of the modes, has the form
$\tens{T}=\sum_{i=1}^{r(\tens{T})}\vect{a}_i\otimes\vect{v}_i^{\otimes2}$.
A decomposition contains no cubes if none of its summands is
proportional to a cube.  In a decomposition of the displayed form, this means that
$\vect{a}_i$ and $\vect{v}_i$ are not proportional for any $i$.

The argument using wedge supports begins with the following bounds for square covers in dimensions two and three. 
The ternary assertions over $\RR$ are
\cite[Lemma~5.4]{li2026comonrankgap}.

\begin{proposition}[Square covers in small dimensions]\label{lem:square-covers}
Let $\set{W}\subseteq\Sym^2(\FF^n)$.
There exist vectors $\vect{\ell}_1,\ldots,\vect{\ell}_N\in\FF^n$
such that
$\set{W}\subseteq\Span\{\vect{\ell}_1^2,\ldots,\vect{\ell}_N^2\},$ 
where the following choices of $N$ suffice:
\[
N=
\begin{cases}
3,&n=2,\\
4,&n=3,\ \dim\set{W}\le2,\\
5,&n=3,\ \dim\set{W}\le4.
\end{cases}
\]
If $\set{W}=\Span\{Q\}$, one may take $N=\rank Q$.
\end{proposition}

\begin{proof}
Congruence diagonalization gives
$Q=\sum_{i=1}^{\rank Q}\lambda_i\vect{\ell}_i^2$,
which proves the assertion for $\dim\set{W}\le1$.
The binary assertion follows from
$\Sym^2(\FF^2)=\Span\{x^2,y^2,(x+y)^2\}$. 
Suppose $n=3$ and $\dim\set{W}=2$.
In matrix coordinates, write
$\set{W}=\Span\{\matr{B},\matr{C}\}$.
If every matrix in $\set{W}$ is singular, diagonalization gives
a cover by at most $\rank\matr{B}+\rank\matr{C}\le4$ squares.
Otherwise, choose $\matr{B}$ nonsingular.
The cubic $\det(t\matr{B}+\matr{C})$ has a root $t_0\in\FF$.
The matrix $\matr{M}=t_0\matr{B}+\matr{C}$ then satisfies
$\set{W}=\Span\{\matr{B},\matr{M}\}$ and
$1\le\rank\matr{M}\le2$.
If $\rank\matr{M}=1$, diagonalization again gives a cover
by four squares.

Suppose $\rank\matr{M}=2$.
After simultaneous congruence, we may assume
$\matr{M}=\mathsf{diag}(m_1,m_2,0)$ with $m_1m_2\ne0$.
On $\set{U}=\Span\{\vect{e}_1,\vect{e}_2\}$, define
$q_1(\vect{u})=u_1^2/m_1+u_2^2/m_2$ and
$q_2(\vect{u})=\vect{u}^{\mathsf T}\matr{B}^{-1}\vect{u}$.
The first quadratic is nonzero.
If the second vanished identically, orthogonality for
$\matr{B}^{-1}$ would give $\set{U}\subseteq\set{U}^{\perp}$,
contradicting $\dim\set{U}=2>\dim\set{U}^{\perp}=1$.
Since $\FF$ is infinite, there is $\vect{u}\in\set{U}$
such that $q_1(\vect{u})q_2(\vect{u})\ne0$.
Set $\matr{R}=\vect{u}\vect{u}^{\mathsf T}$.
Apply the matrix determinant lemma to the nonsingular leading
$2\times2$ block of $\matr{M}$ and to $\matr{B}$.
Since a rank-one perturbation changes rank by at most one, this gives
\begin{align*}
\rank\bigl(\matr{M}-q_1(\vect{u})^{-1}\matr{R}\bigr)=1,\ \
\rank\bigl(\matr{B}-q_2(\vect{u})^{-1}\matr{R}\bigr)=2.
\end{align*}
These two matrices and $\matr{R}$ span a space containing
$\set{W}$, so diagonalization gives a cover by
$1+2+1=4$ squares.

Now suppose $\set{W}\subseteq\Sym^2(\FF^3)$ and
$\dim\set{W}\le4$.
Its annihilator $\set{A}$ in the dual quadratic space has
dimension $6-\dim\set{W}\ge2$.
Over $\CC$, $\set{A}$ contains a quadratic of rank at least two.
Otherwise, any two linearly independent nonzero elements of $\set{A}$ would have rank one, whereas their sum would have rank two, a contradiction. 
Over $\RR$, suppose $\set{A}$ contains no indefinite quadratic.
After changing signs, choose independent positive semidefinite
forms $P,R\in\set{A}$.
The quotient $R/P$ is nonconstant on $\{P>0\}$, since constancy
would imply proportionality.
Thus there are $\vect{x},\vect{y}$ with
$P(\vect{x}),P(\vect{y})>0$ and a scalar $\lambda$ satisfying
$R(\vect{x})/P(\vect{x})<\lambda<R(\vect{y})/P(\vect{y})$.
This gives
$(R-\lambda P)(\vect{x})<0<(R-\lambda P)(\vect{y})$,
a contradiction.

Over either field, congruence and scalar multiplication reduce
a suitable element of $\set{A}$ to
$Q_1=x^2-y^2$ or $Q_2=x^2+y^2-z^2$.
Consider the corresponding sets
\begin{align*}
\set{N}_1
&=\{(0,0,1),(1,1,0),(1,-1,0),(1,1,1),(1,-1,1)\},\\
\set{N}_2
&=\{(1,0,1),(-1,0,1),(0,1,1),(0,-1,1),(3,4,5)\}.
\end{align*}
For each $j\in\{1,2\}$, direct calculation gives
$Q_j(\vect{v})=0$ for every $\vect{v}\in\set{N}_j$ and
$\dim\Span\{\vect{v}^2:\vect{v}\in\set{N}_j\}=5$.
Consequently, the transformed space satisfies
$\set{W}\subseteq\{Q_j\}^{\perp}
=\Span\{\vect{v}^2:\vect{v}\in\set{N}_j\}$.
Transforming back gives the required cover by five squares.
\end{proof}

A symmetric tensor on $\set{V}$ is called concise if its essential
space is $\set{V}$.
For a bivector in $\bigwedge^2\set{V}$, its rank and support are the
rank and image of the contraction map $\set{V}^*\to\set{V}$, respectively.
Let $\tens{S}\in\Sym^3(\set{V})$ be concise and choose a
minimum CP decomposition
\[
\tens{S}=\sum_{i=1}^{r(\tens{S})}
\vect{a}_i\otimes\vect{b}_i\otimes\vect{c}_i.
\]
Define $\kappa_{ab}\eqdef\dim\Span\{\vect{a}_i\wedge\vect{b}_i\}$,
let $\set{H}_{ab}$ be the sum of the supports of its nonzero wedges,
and put $h_{ab}\eqdef\dim\set{H}_{ab}$.  Define the other two mode pairs cyclically.
These quantities refer to the chosen decomposition.
Symmetry gives
\begin{equation}\label{eq:wedge-identities}
\sum_i(\vect{a}_i\wedge\vect{b}_i)\otimes\vect{c}_i=0.
\end{equation}
The matrix of the remaining factors has rank $d(\tens{S})$; its row space is contained
in the kernel of the wedge map.  A basis chosen from the displayed
wedges spans the wedge space.  Therefore, the wedge dimensions satisfy
\begin{equation}\label{eq:wedge-budget}
\kappa_{ab},\kappa_{bc},\kappa_{ca}\le e_r(\tens{S}),\ \ 
h_{ab}\le2\kappa_{ab}.
\end{equation}
This is the wedge dimension estimate of
\cite[Lemma~4.2]{li2026comonrankgap}.

We will also use the following support-dimension estimate.  If
$\tens{\omega}=\sum_{i=1}^m\lambda_i
(\vect{a}_i\wedge\vect{b}_i)$ has rank $2r$ and every
$\lambda_i\ne0$, write its matrix as
$\matr{\Omega}=\matr{X}\matr{J}\matr{X}^{\mathsf T}$, using the convention
$\vect{a}\wedge\vect{b}\leftrightarrow\vect{a}\vect{b}^{\mathsf T}-\vect{b}\vect{a}^{\mathsf T}$.
Here $\matr{X}=[\vect{a}_1\ \vect{b}_1\ \cdots\ \vect{a}_m\ \vect{b}_m]$, and
$\matr{J}$ is the block-diagonal alternating $2m\times2m$ matrix with
$2\times2$ blocks and $(2i-1,2i)$ entry $\lambda_i$; it is nonsingular.
The radical of $\mathsf{row}\matr{X}$ with respect to the alternating form defined by $\matr{J}$ has dimension
$\rank\matr{X}-2r\le2m-\rank\matr{X}$.
Thus the dimension of the span of all factor vectors,
$\rank\matr{X}=\dim\Span\{\vect{a}_i,\vect{b}_i:1\le i\le m\}$, satisfies 
\begin{equation}\label{eq:endpoint-rank}
\rank\matr{X}\le m+r.
\end{equation}

A circuit is a minimally linearly dependent subfamily.
The following ordering allows us to estimate supports one circuit
at a time.

\begin{proposition}[An ordering by circuits]\label{lem:circuit-ear}
Let $(\vect{x}_e)_{e\in \set{E}}$ be a finite family of nonzero vectors whose
linear matroid is connected, and fix $e_0\in \set{E}$.  There are sets
$\{e_0\}=\set{A}_0\subsetneq \set{A}_1\subsetneq\cdots\subsetneq \set{A}_t=\set{E}$ 
and circuits $\set{C}_1,\ldots,\set{C}_t$ such that, with
$\set{N}_j=\set{C}_j\setminus \set{A}_{j-1}$, the following relations hold:
\[
\set{A}_j=\set{A}_{j-1}\cup\set{C}_j,
\ 
\set{C}_j\cap\set{A}_{j-1}\ne\varnothing,
\ 
\set{N}_j\ne\varnothing,
\]
and
\begin{equation}\label{eq:circuit-ear-rank}
\dim
\frac{\Span\{\vect{x}_e:e\in \set{A}_j\}}
     {\Span\{\vect{x}_e:e\in \set{A}_{j-1}\}}
=|\set{N}_j|-1.
\end{equation}
\end{proposition}

\begin{proof}
Suppose that a nonempty proper set $\set{A}\subset \set{E}$ has already been
constructed.  Connectedness gives a circuit meeting both $\set{A}$ and
$\set{E}\setminus \set{A}$.  Choose such a circuit $\set{C}$ for which
$|\set{C}\setminus \set{A}|$ is minimal, and put $\set{N}=\set{C}\setminus \set{A}$. 
The circuit relation, whose coefficients are all nonzero, gives a
nontrivial relation among the images of
$\{\vect{x}_e:e\in \set{N}\}$ modulo $\Span\{\vect{x}_e:e\in \set{A}\}$.  Thus their
quotient rank is at most $|\set{N}|-1$.  If it were at most $|\set{N}|-2$, the quotient
relation space would have dimension at least two.  A linear combination
of an independent quotient relation and the circuit relation would then
give a nonzero quotient relation supported on a proper subset
$\set{N}'\subsetneq \set{N}$.  Since $\set{N}'$ is a proper subset of the circuit $\set{C}$,
its vectors are independent; their nonzero combination in the span of
$\set{A}$ is therefore not the zero vector.  Choose a basis of
$\Span\{\vect{x}_e:e\in \set{A}\}$ from the vectors indexed by $\set{A}$.  That basis,
together with the vectors indexed by $\set{N}'$, contains a circuit meeting
both parts, contrary to the minimality of $|\set{N}|$.  Therefore, the quotient
rank is $|\set{N}|-1$.  Replace $\set{A}$ by $\set{A}\cup \set{C}$ and repeat.  Since $\set{E}$ is finite, the
process ends after finitely many steps.
\end{proof}

This ordering bounds the total support of connected wedge families.

\begin{proposition}[Connected wedge supports]\label{lem:connected-wedge}
Let a finite family of nonzero decomposable bivectors span a space
of dimension $\kappa$.  If its linear matroid is connected, its total
support dimension is at most $\kappa+1$.
\end{proposition}

\begin{proof}
Write the bivectors as
$\tens{\omega}_e=\vect{a}_e\wedge\vect{b}_e$ and apply
\cref{lem:circuit-ear}.  For $\set{A}\subseteq \set{E}$, set
$\set{W}(\set{A})\eqdef\sum_{e\in \set{A}}\Span\{\vect{a}_e,\vect{b}_e\}$ and
$\set{L}(\set{A})\eqdef\Span\{\tens{\omega}_e:e\in \set{A}\}.$ 
Thus $\dim\set{W}(\set{A}_0)=2$ and $\dim\set{L}(\set{A}_0)=1$. 
At step $j$, the circuit relation gives
\[
\tens{\eta}_j
=\sum_{e\in \set{N}_j}\alpha_e\tens{\omega}_e
=-\sum_{e\in \set{C}_j\cap \set{A}_{j-1}}\alpha_e\tens{\omega}_e
\in\bigwedge^2\set{W}(\set{A}_{j-1}),
\]
where every $\alpha_e$ is nonzero.  
Since $\set{C}_j$ meets $\set{A}_{j-1}$,
$\set{N}_j$ is a nonempty proper subset of $\set{C}_j$.
The corresponding bivectors are therefore linearly independent.
As every $\alpha_e$ is nonzero, it follows that
$\tens{\eta}_j\ne0$.   
Suppose that $\tens{\eta}_j$ has
rank $2r$, where $r\ge1$, and let $h_j$ be the total support dimension
of the new bivectors.  The support-dimension estimate
\eqref{eq:endpoint-rank} gives $h_j\le|\set{N}_j|+r$.  The $2r$-dimensional
support of $\tens{\eta}_j$ is contained both in $\set{W}(\set{A}_{j-1})$ and in the
support of the new batch.  Consequently, the support increment satisfies
\[
\dim\set{W}(\set{A}_j)-\dim\set{W}(\set{A}_{j-1})
\le h_j-2r
\le|\set{N}_j|-r
\le|\set{N}_j|-1.
\]
By \eqref{eq:circuit-ear-rank}, the corresponding increase of
$\dim\set{L}(\set{A}_j)$ is exactly $|\set{N}_j|-1$.  Induction gives
$$\dim\set{W}(\set{E})-2\le\dim\set{L}(\set{E})-1,$$ and hence
$\dim\set{W}(\set{E})\le\kappa+1$.
\end{proof}

We next combine the support bounds with square covers to remove small
wedge components without increasing the decomposition length.

\begin{proposition}[Elimination of small wedge components]
\label{lem:small-components}
Let $\tens{T}\in\Sym^2(\set{V})\otimes\set{A}$.
Among its minimum CP decompositions
$\sum_i\vect{a}_i\otimes\vect{b}_i\otimes\vect{c}_i$,
choose one minimizing $\kappa_{ab}$.
The nonzero $ab$-wedge matroid then has no connected component of
matroid rank one or two.  Consequently, its support satisfies
\begin{equation}\label{eq:small-component-support}
h_{ab}\le\left\lfloor\frac{4\kappa_{ab}}3\right\rfloor.
\end{equation}
\end{proposition}
\begin{proof}
Set $\tens{\omega}_i\eqdef\vect{a}_i\wedge\vect{b}_i$, and let $\set{E}_\nu$ be a
connected component of the linear matroid of the nonzero $\tens{\omega}_i$.
Write
$\set{W}_\nu\eqdef\Span\{\tens{\omega}_i:i\in \set{E}_\nu\}$, 
$\kappa_\nu\eqdef\dim\set{W}_\nu$ and 
$m_\nu\eqdef|\set{E}_\nu|.$ 
The spaces $\set{W}_\nu$ form a direct sum.
For each $\nu$, choose a basis of $\set{W}_\nu$ from the wedges
indexed by $\set{E}_\nu$.
If the union of these bases were linearly dependent, it would
contain a circuit meeting at least two components, a contradiction.  
Projecting \eqref{eq:wedge-identities} onto
$\set{W}_\nu$ gives
\begin{equation}\label{eq:component-wedge-identity}
\sum_{i\in \set{E}_\nu}\tens{\omega}_i\otimes\vect{c}_i=0.
\end{equation}
Choose a basis of $\set{W}_\nu$ and write the $\tens{\omega}_i$ in that basis.
The resulting $\kappa_\nu\times m_\nu$ coefficient matrix has row rank
$\kappa_\nu$.  Therefore, \eqref{eq:component-wedge-identity} gives
$\kappa_\nu$ independent relations among the vectors $\vect{c}_i$, and
\begin{equation*}
\dim\Span\{\vect{c}_i:i\in \set{E}_\nu\}
\le m_\nu-\kappa_\nu.
\end{equation*}
Since every $\vect{c}_i$ is nonzero, this also implies
$m_\nu\ge\kappa_\nu+1$.
Define the component subtotal by
$\tens{T}_\nu
\eqdef\sum_{i\in \set{E}_\nu}
 \vect{a}_i\otimes\vect{b}_i\otimes\vect{c}_i.$
Let $\set{L}_\nu$ be the total support of its wedges.  Equation \eqref{eq:component-wedge-identity} says that
$\tens{T}_\nu$ is symmetric in its first two modes.  Its quadratic slice
space $\set{Q}_\nu\subseteq\Sym^2(\set{L}_\nu)$ consequently satisfies
$\dim\set{Q}_\nu\le m_\nu-\kappa_\nu$.

Suppose first that $\kappa_\nu=1$.  Then
$\dim\set{L}_\nu\le2$ by \cref{lem:connected-wedge}, and $m_\nu\ge2$.
If $m_\nu=2$, the slice space has dimension at most one, and its single
binary quadratic form is a linear combination of at most two squares of linear forms.  If $m_\nu\ge3$, use
the binary cover by three squares in \cref{lem:square-covers}.  Thus
$\tens{T}_\nu$ has a CP decomposition that is symmetric in the first two
modes and has at
most $m_\nu$ terms.

Suppose next that $\kappa_\nu=2$.  The bounds are then
$\dim\set{L}_\nu\le3,\ m_\nu\ge3$ and 
$\dim\set{Q}_\nu\le m_\nu-2.$ 
For $m_\nu=3$, a single ternary
quadratic form is a linear combination of at most three squares of linear forms.  For $m_\nu=4$, use the
cover by four squares of a ternary pencil.  For $m_\nu=5$, use the
cover by five squares in \cref{lem:square-covers}.  Finally, for
$m_\nu\ge6$, use a basis of six squares of $\Sym^2(\FF^3)$.  Hence, in
every case, there is an expression
\[
\tens{T}_\nu
=\sum_{j=1}^{t_\nu}
(\vect{v}_{\nu j})^2\otimes\vect{d}_{\nu j},
\ \ 
t_\nu\le m_\nu.
\]

If $t_\nu<m_\nu$, replacing the component subtotal gives a CP decomposition shorter than the chosen minimum decomposition.  Thus $t_\nu=m_\nu$.  
The resulting full decomposition has the same length and is still a minimum CP decomposition, while all $ab$-wedges of the replacement terms vanish.  
Since $\set{W}_\nu$ is a direct summand of
the original wedge span, the replacement decreases $\kappa_{ab}$ by
$\kappa_\nu$, contrary to the choice of the decomposition.  Therefore,
no component with $\kappa_\nu=1$ or $2$ occurs.

All remaining components have $\kappa_\nu\ge3$.  By
\cref{lem:connected-wedge}, their support dimensions satisfy
$\dim\set{L}_\nu\le\kappa_\nu+1$.  Consequently, their total contribution obeys
\[
h_{ab}\le\sum_\nu\dim\set{L}_\nu
\le\sum_\nu(\kappa_\nu+1)
\le\sum_\nu\frac{4\kappa_\nu}{3}
=\frac{4\kappa_{ab}}3.
\]
Since $h_{ab}$ is an integer, this proves
\eqref{eq:small-component-support}.
\end{proof}

For a decomposition containing no cubes, the three mode pairs also give
a dimension count and a symmetric rank bound. The count of terms
with three nonzero wedges follows
\cite[Lemma~4.8]{li2026comonrankgap}.

\begin{proposition}[A count using the three mode pairs]\label{lem:three-sector}
Let $\tens{S}\in\Sym^3(\set{V})$ be concise and have a minimum CP decomposition
containing no cubes.  Then
\begin{align}
2d(\tens{S})\le h_{ab}+h_{bc}+h_{ca},
\ \
s(\tens{S})\le d(\tens{S})+4e_r(\tens{S})+h_{ab}+h_{bc}+h_{ca}.
\label{eq:three-sector-rank}
\end{align}
\end{proposition}
\begin{proof}
Each factor lies in the intersection of two spaces of wedge supports.  For
example, if $\vect{a}_i\wedge\vect{b}_i$ and
$\vect{a}_i\wedge\vect{c}_i$ are nonzero, then
$\vect{a}_i\in\set{H}_{ab}\cap\set{H}_{ca}$.  If one wedge vanishes,
the absence of a cube makes the other two nonzero and gives the same conclusion
for the proportional factors.  Hence the three pairwise intersections
span the actual support.  If their dimensions are
$t_{12},t_{23},t_{31}$ and their common intersection has dimension
$t$, then
\[
d(\tens{S})\le t_{12}+t_{23}+t_{31}-2t,
\ \ 
\sum h_{xy}\ge2(t_{12}+t_{23}+t_{31})-3t\ge2d(\tens{S}).
\]

Partition the terms into the three classes having one vanishing wedge
and the class having no vanishing wedge; let the last class have
cardinality $x$.  Modulo $\set{H}_{bc}$, the vectors $\vect{b}_i$ for
which $\vect{b}_i\wedge\vect{c}_i=0$ span the quotient.  Thus there are at least
$d(\tens{S})-h_{bc}$ such indices, and cyclically the same holds for the other
two classes.  Since the total number of terms is $r(\tens{S})$, this count gives
\[
x\le e_r(\tens{S})-2d(\tens{S})+h_{ab}+h_{bc}+h_{ca}.
\]
The polarization identities in \cref{sec:quaternary-rank-bound}
express a product with a repeated factor as a sum of at most three
cubes, and a product of three arbitrary factors as a sum of at most
four cubes.  Symmetrizing the decomposition therefore gives
$s(\tens{S})\le3r(\tens{S})+x$, which is
\eqref{eq:three-sector-rank}.
\end{proof}

\subsubsection{Proof of rank equality}
\begin{proof}[Proof of \cref{prop:e3}]
For $\FF=\RR$, the assertion under $e_r(\tens{S})\le3$ is
\cite[Proposition~3.3]{li2026comonrankgap}.
For $\FF=\CC$, we use the proof of that proposition with the
following observations. The square cover bounds needed in
\cite[Lemma~5.5]{li2026comonrankgap} hold over $\CC$ by
\cref{lem:square-covers}, and quadratic forms admit diagonalization
by congruence. The bounds for symmetric rank in dimensions at most
four used in
\cite[Propositions~5.6 and~5.9]{li2026comonrankgap}
hold over $\CC$ by \eqref{eq:rank-cap}.
The decomposable-bivector support and direction arguments in
\cite[Lemma~5.7]{li2026comonrankgap} use exterior products,
linear dependence, and dimension counts. The polar-block splitting in
\cite[Lemma~5.8]{li2026comonrankgap} follows from commutativity
of contractions and the associated dual projections.
Both arguments remain valid over $\CC$.
In the proof of \cite[Proposition~5.9]{li2026comonrankgap},
the polar space is identified with the kernel of the wedge map;
the subsequent argument either produces a transverse rank-one polar
or splits the tensor into binary and ternary blocks.
These steps use the preceding support and splitting arguments,
while the rank bounds for the resulting blocks follow from
\eqref{eq:rank-cap}.
Together with \cref{lem:square-covers} and congruence diagonalization,
this gives the complex versions of
\cite[Lemma~5.5 and Propositions~5.6 and~5.9]{li2026comonrankgap}.

We can therefore repeat the induction on CP rank in
\cite[proof of Proposition~3.3]{li2026comonrankgap}.
Its reductions either give rank equality directly or express
$\tens{S}$ as a sum of $q\ge1$ cubes and a symmetric tensor $\tens{S}'$
with $r(\tens{S}')=r(\tens{S})-q$ and
$e_r(\tens{S}')\le e_r(\tens{S})\le3$.
For removal of displayed cubes these relations follow from
\cite[Lemma~2.4]{li2026comonrankgap}; for removal of a cube obtained
from a transverse polar they follow from
\cite[Lemma~2.5 and Proposition~5.9]{li2026comonrankgap}.
Both removal arguments use only minimality, projections, and
dimension counts, and every complex scalar has a cube root.
The zero tensor starts the induction, and the induction hypothesis
gives 
$$s(\tens{S})\le q+s(\tens{S}')=q+r(\tens{S}')=r(\tens{S}).$$
The reverse inequality always holds. 
Finally, if $e_s(\tens{S})\le4$ and $\Delta_c(\tens{S})>0$,
integrality gives
$e_r(\tens{S})=e_s(\tens{S})-\Delta_c(\tens{S})\le3$.
The rank equality just proved yields a contradiction.
\end{proof}

% \begin{remark}[Relation to known ranges of rank equality]
% Friedland's flattening criterion gives rank equality when $e_r\le1$
% \cite[Theorem~1.1]{friedland2016remarks}; his argument for removing cubes
% when the CP rank excess equals one is given in
% \cite[Lemma~5.1]{friedland2016remarks}.
% Seigal proves rank equality over $\CC$ for cubics in at most four
% variables \cite[Theorem~1.3]{seigal2020ranks}.
% The real range $e_r\le3$ is established in
% \cite[Proposition~3.3]{li2026comonrankgap}.
% The extension described above gives the same range over $\CC$
% without a dimension restriction.
% \end{remark}

\section{Proofs of the symmetric rank gap bounds}
\label{sec:gap-proofs}
We first prove the bounds based on symmetric rank excess and the
Comon rank gap.  We then use the CP rank excess estimates from
\cref{sec:rank-inputs} to control the symmetric rank gap by CP rank excess.
Finally, polar extraction gives the dimension bound, and the symmetric
rank excess estimate gives additivity in low dimensions.

\subsection{Proof of the one-sided bounds}
\label{sec:symmetric-elementary}

We begin by adjoining cubes in a complementary variable space.

\begin{proposition}[Adjoining cubes of linearly independent linear forms]\label{lem:independent-cube}
Let $\tens{S}$ be a symmetric third-order tensor and let
$\vect{y}_1,\ldots,\vect{y}_h$ be linearly independent vectors such that
$\Span\{\vect{y}_1,\ldots,\vect{y}_h\}\cap\Ess(\tens{S})=\{0\}$.
Then
\[
s\left(\sum_{j=1}^h\vect{y}_j^3\oplus\tens{S}\right)
=h+s(\tens{S}).
\]
\end{proposition}
\begin{proof}
It suffices to prove the assertion for one cube and then argue
by induction.
The one-cube assertion is immediate when $\tens{S}=0$,
so assume $\tens{S}\ne0$.
Projection and subadditivity give
$s(\tens{S})\le s(\vect{y}^3\oplus\tens{S})
\le s(\tens{S})+1$.
Suppose that equality holds in the first inequality and write 
$$\vect{y}^3\oplus\tens{S}
=\sum_{i=1}^{s(\tens{S})}(\alpha_i\vect{y}+\vect{\ell}_i)^3.$$  
Here each $\vect{\ell}_i$ belongs to $\Ess(\tens{S})$ by
\cref{lem:minimum}.
If $\alpha_i\notin\{0,1\}$ for some $i$, apply the linear map that is
the identity on $\Ess(\tens{S})$ and sends
$\vect{y}$ to $\vect{\ell}_i/(1-\alpha_i)$.  The $i$th cube then equals the cube
on the left, and cancellation expresses $\tens{S}$ with at most
$s(\tens{S})-1$ cubes.
If every $\alpha_i$ is zero or one, comparison of the
$\vect{y}^3$ coefficients shows that precisely one coefficient is one.
Comparison of the $\vect{y}^2$ coefficients then gives
$\vect{\ell}_i=0$ for that index, with the same contradiction.
This proves the assertion for one cube.  Repeated application gives the
result for $h$ cubes.
\end{proof}

% \begin{remark}[Relation to additivity for cubes of linearly independent linear forms]
% Over $\CC$, the case of one cube is contained in
% \cite[Proposition~3.1 and Corollary~3.2]{carlini2015symmetricstrassen}.
% The same substitution argument works for real cubics, and iteration
% gives the form used here.
% \end{remark}

We next extract cubes with linearly independent factors from one summand to bound the gap.

\begin{proposition}[Elementary one-sided bounds]
\label{lem:elementary-one-sided}
For nonzero symmetric tensors $\tens{S}_1$ and $\tens{S}_2$, one has
\begin{align}
0\le\Delta_s(\tens{S}_1,\tens{S}_2)
\le2\min\{e_s(\tens{S}_1),e_s(\tens{S}_2)\},
\label{eq:elementary-excess}\ \
\Delta_s(\tens{S}_1,\tens{S}_2)
\le\min\{s(\tens{S}_1),s(\tens{S}_2)\}-1. 
\end{align}
\end{proposition}

\begin{proof}
Choose $d(\tens{S}_1)$ independent factors in a minimum Waring
decomposition and let $\tens{S}_{1,\mathrm{ind}}$ be their subtotal.  Write
$\tens{S}_{1,\mathrm{rem}}\eqdef\tens{S}_1-\tens{S}_{1,\mathrm{ind}}$.  Then
$s(\tens{S}_{1,\mathrm{ind}})=d(\tens{S}_1)$ and 
$s(\tens{S}_{1,\mathrm{rem}})\le e_s(\tens{S}_1).$ 
By \cref{lem:independent-cube}, the symmetric ranks satisfy
\begin{align*}
d(\tens{S}_1)+s(\tens{S}_2)
=s(\tens{S}_{1,\mathrm{ind}}\oplus\tens{S}_2) 
\le s(\tens{S}_1\oplus\tens{S}_2)+s(-\tens{S}_{1,\mathrm{rem}}) 
\le s(\tens{S}_1\oplus\tens{S}_2)+e_s(\tens{S}_1).
\end{align*}
Using $s(\tens{S}_1)=d(\tens{S}_1)+e_s(\tens{S}_1)$ gives
$\Delta_s(\tens{S}_1,\tens{S}_2)\le2e_s(\tens{S}_1)$. 
For the symmetric rank bound, choose $\alpha\in\Ess(\tens{S}_1)^*$ such that
$\tens{S}_1(\alpha,\alpha,\alpha)\ne0$, and let $\vect{z}$ span a new one-dimensional space.  
Extend the map
$\vect{x}\mapsto\alpha(\vect{x})\vect{z}$ on $\Ess(\tens{S}_1)$
by the identity on $\Ess(\tens{S}_2)$.
The induced map on symmetric tensors sends $\tens{S}_1$ to the nonzero
cube $\tens{S}_1(\alpha,\alpha,\alpha)\vect{z}^3$
and leaves $\tens{S}_2$ unchanged. 
Therefore,
one has
$$s(\tens{S}_1\oplus\tens{S}_2)
\ge s\bigl(\tens{S}_1(\alpha,\alpha,\alpha)\vect{z}^3\oplus\tens{S}_2\bigr)
=1+s(\tens{S}_2).$$  
Interchanging the summands completes both estimates.
\end{proof}

The next proposition replaces one summand by a tensor supported on
the projected mixed factors while preserving a positive gap.

\begin{proposition}[Replacing a summand by its mixed part]\label{lem:projected-core}
If $\Delta_s(\tens{S}_1,\tens{S}_2)>0$, there is a nonzero symmetric
tensor $\widehat{\tens{S}}_1$ such that
\begin{align}\label{eq:projected-core}
d(\widehat{\tens{S}}_1)\le e_s(\tens{S}_1)+\Delta_c(\tens{S}_2)
-\Delta_s(\tens{S}_1,\tens{S}_2),\ \
\Delta_s(\widehat{\tens{S}}_1,\tens{S}_2)\ge\Delta_s(\tens{S}_1,\tens{S}_2). 
\end{align}
\end{proposition}
\begin{proof}
Identify the direct sum with a tensor on
$\set{A}\oplus\set{B}$, where
$\set{A}\eqdef\Ess(\tens{S}_1)$ and $\set{B}\eqdef\Ess(\tens{S}_2)$, and choose a minimum Waring decomposition 
\[
\tens{S}_1\oplus\tens{S}_2
=\sum_{i=1}^t(\vect{x}_i+\vect{y}_i)^3,
\ \ 
\vect{x}_i\in\set{A},
\ 
\vect{y}_i\in\set{B}.
\]
Set
$\set{E}\eqdef\Span\{\vect{x}_i:\vect{y}_i\ne0\}$ and 
$\set{M}\eqdef\Span\{(\vect{x}_i+\vect{y}_i)\otimes\vect{y}_i:
\vect{y}_i\ne0\}.$ 
The slice space of $\tens{S}_2$ lies in $\set{M}$, so the support in the
first factor of $\set{M}$ contains $\set{B}$.
It also contains $\vect{x}_i+\vect{y}_i$ whenever $\vect{y}_i\ne0$;
subtracting $\vect{y}_i\in\set{B}$ shows that it contains $\set{E}$.
The reverse inclusion follows from the generators of $\set{M}$,
so this support is $\set{B}\oplus\set{E}$.
Therefore, \cref{lem:ambient-support} gives
$\dim\set{M}\ge r(\tens{S}_2)+\dim\set{E}$.  On
$\set{L}\eqdef\Span\{(\vect{x}_i+\vect{y}_i)^2\}$, consider the map
$\pi\eqdef(\id\otimes\mathsf{pr}_{\set{B}})|_{\set{L}}.$ 
Its image is $\set{M}$, and its kernel contains the slice space of
$\tens{S}_1$.  Hence
\begin{align*}
d(\tens{S}_1)+r(\tens{S}_2)+\dim\set{E}\le\dim\set{L}\le t,\ \
\dim\set{E}\le e_s(\tens{S}_1)+\Delta_c(\tens{S}_2)
-\Delta_s(\tens{S}_1,\tens{S}_2).
\end{align*}
Set $b\eqdef\#\{i:\vect{y}_i=0\}$ and
$\widehat{\tens{S}}_1\eqdef\sum_{\vect{y}_i\ne0}\vect{x}_i^3$.
The mixed terms cancel, and \cref{lem:minimum} gives
$s(\widehat{\tens{S}}_1\oplus\tens{S}_2)=t-b$.  Moreover, the symmetric rank satisfies
$s(\widehat{\tens{S}}_1)\ge s(\tens{S}_1)-b$.  It follows that
$$\Delta_s(\widehat{\tens{S}}_1,\tens{S}_2)\ge\Delta_s(\tens{S}_1,\tens{S}_2)>0$$ and
$d(\widehat{\tens{S}}_1)\le\dim\set{E}$.  In particular, the tensor $\widehat{\tens{S}}_1$
is nonzero.
\end{proof}

A symmetric rank bound in small dimensions now turns this reduction into a
descent argument.

\begin{proposition}[Descent from a symmetric rank bound in small dimensions]\label{lem:low-dimensional-descent}
Let $h\ge1$ be an integer and suppose 
\begin{equation}\label{eq:descent-hypothesis}
M_{3,\FF}(n)\le2n-1\ \ (1\le n\le h).
\end{equation}
Then
\begin{equation}\label{eq:general-descent}
\Delta_s(\tens{S}_1,\tens{S}_2)\le
\begin{cases}
(\Delta_c(\tens{S}_2)-1)_+,&e_s(\tens{S}_1)\le h,\\
(\Delta_c(\tens{S}_2)+e_s(\tens{S}_1)-h-1)_+,&e_s(\tens{S}_1)\ge h+1.
\end{cases}
\end{equation}
\end{proposition}
\begin{proof}
The assertion is immediate when the gap is zero.  Otherwise, suppose,
for contradiction, that
$\Delta_s(\tens{S}_1,\tens{S}_2)
\ge\Delta_c(\tens{S}_2)+(e_s(\tens{S}_1)-h)_+.$ 
By \cref{lem:projected-core}, there is $\widehat{\tens{S}}_{1,0}\ne0$ with
\begin{align*}
d(\widehat{\tens{S}}_{1,0})\le e_s(\tens{S}_1)+\Delta_c(\tens{S}_2)
-\Delta_s(\tens{S}_1,\tens{S}_2)\le h,\ \
\Delta_s(\widehat{\tens{S}}_{1,0},\tens{S}_2)\ge\Delta_s(\tens{S}_1,\tens{S}_2)
\ge\Delta_c(\tens{S}_2).
\end{align*}
Given $\widehat{\tens{S}}_{1,i}\ne0$ with $d(\widehat{\tens{S}}_{1,i})\le h$,
equations~\eqref{eq:descent-hypothesis} and
\eqref{eq:projected-core} give
\begin{align*}
e_s(\widehat{\tens{S}}_{1,i})&\le d(\widehat{\tens{S}}_{1,i})-1,\ \
d(\widehat{\tens{S}}_{1,i+1})
\le e_s(\widehat{\tens{S}}_{1,i})+\Delta_c(\tens{S}_2)
-\Delta_s(\widehat{\tens{S}}_{1,i},\tens{S}_2)
\le d(\widehat{\tens{S}}_{1,i})-1,\\
\Delta_s(\widehat{\tens{S}}_{1,i+1},\tens{S}_2)
&\ge\Delta_s(\widehat{\tens{S}}_{1,i},\tens{S}_2) 
\ge\Delta_s(\tens{S}_1,\tens{S}_2)>0.
\end{align*}
This gives an impossible strictly decreasing sequence
$d(\widehat{\tens{S}}_{1,0})>d(\widehat{\tens{S}}_{1,1})>\cdots\ge1$.
Therefore, every positive gap satisfies
$\Delta_s(\tens{S}_1,\tens{S}_2)
<\Delta_c(\tens{S}_2)+(e_s(\tens{S}_1)-h)_+.$ 
Since all quantities are integers, this is
\eqref{eq:general-descent}.
\end{proof}

\begin{proof}[Proof of \cref{prop:one-sided}]
\Cref{lem:elementary-one-sided} gives the lower bound, the term
$2e_s(\tens{S}_1)$ in both branches of
\eqref{eq:one-sided-piecewise}, and
\eqref{eq:minimum-rank-bound}.  Apply
\cref{lem:low-dimensional-descent} with $h=4$ to
$(\tens{S}_1,\tens{S}_2)$, using
$M_{3,\FF}(n)\le2n-1$ for $1\le n\le4$ from
\eqref{eq:rank-cap}.  This gives the entries involving
$\Delta_c(\tens{S}_2)$.  If $r(\tens{S}_2)=s(\tens{S}_2)$, then
$\Delta_c(\tens{S}_2)=0$, and the same estimate gives
\eqref{eq:eq-partner}.
\end{proof}

% \begin{remark}[Relation to symmetric additivity criteria]
% Carlini, Catalisano, and Chiantini prove additivity for a pure power
% summand and pairs of binary forms over $\CC$
% \cite[Proposition~3.1 and Theorem~5.6]{carlini2015symmetricstrassen}.
% Here the projected mixed contribution gives a dimension reduction
% that also applies when the gap is positive.  Combined with the
% symmetric rank bounds in small dimensions, it yields the estimate involving the Comon rank gap of the other summand
% over both $\RR$ and $\CC$.
% \end{remark}

\subsection{Comon rank gap and joint excess estimates}

We next use \cref{prop:e3} to bound the Comon rank gap of a single
tensor.  Substitution into the one-sided estimate then completes
the symmetric rank excess bound.

\begin{proof}[Proof of \cref{lem:comon-gap-bound}]
The assertion is immediate for $\tens{S}=0$.  If
$\Delta_c(\tens{S})>0$, then \cref{prop:e3} gives
$e_r(\tens{S})\ge4$.  Therefore, the identity
$\Delta_c(\tens{S})=e_s(\tens{S})-e_r(\tens{S})$ 
implies \eqref{eq:comon-gap-excess-bound}. 
For $q\in\{9,10\}$, suppose that $s(\tens{S})\le q$.
We claim that 
\begin{equation}\label{eq:rank-nine-ten-excess}
e_s(\tens{S})\le q-5.
\end{equation}
Indeed, this follows from $e_s=s-d$ when $d(\tens{S})\ge5$.
For $1\le d(\tens{S})\le4$, equation \eqref{eq:rank-cap} gives
$e_s(\tens{S})\le d(\tens{S})-1\le3$.
Consequently, \eqref{eq:comon-gap-excess-bound} gives
\eqref{eq:rank-nine-comon-zero} and the additional estimate
\begin{equation}\label{eq:rank-ten-comon-one}
s(\tens{S})\le10
\ \Longrightarrow\ 
\Delta_c(\tens{S})\le1.
\end{equation}
\end{proof}

% \begin{remark}[From rank equality to a Comon rank gap bound]
% Criteria for equality of CP rank and symmetric rank such as
% \cite[Theorem~1.1]{friedland2016remarks} identify cases in which the
% Comon rank gap vanishes.  The argument above applies the
% contrapositive of \cref{prop:e3} to a positive gap.  It gives a
% bound for one symmetric tensor that can be used directly in the
% direct sum bound involving the Comon rank gap of the other summand.
% \end{remark}

\begin{proof}[Proof of \cref{prop:joint-excess}]
Under \eqref{eq:symmetric-excess-order}, one has
$e_s(\tens{S}_1)\le e_s(\tens{S}_2)$.
Substituting \eqref{eq:comon-gap-excess-bound} into
\cref{prop:one-sided} gives $(e_s(\tens{S}_2)-5)_+$ when
$e_s(\tens{S}_1)\le4$ and $e_s(\tens{S}_1)+e_s(\tens{S}_2)-9$ when
$e_s(\tens{S}_1)\ge5$.  These are exactly the two cases in
\eqref{eq:two-sided-excess}.
\end{proof}

\begin{proof}[Proof of \cref{thm:sym-q}]
Take the minimum of the one-sided bounds in \cref{prop:one-sided} and
the joint bound in \cref{prop:joint-excess}.  Under
\eqref{eq:symmetric-excess-order}, the joint bound is
$(e_s(\tens{S}_2)-5)_+$ for $e_s(\tens{S}_1)\le4$ and
$e_s(\tens{S}_1)+e_s(\tens{S}_2)-9$ otherwise.  This proves
\eqref{eq:sym-q-main}. 
Item~\ref{item:symmetric-excess-zero} follows from the elementary
one-sided bound, and item~\ref{item:symmetric-excess-rectangle}
follows from \cref{prop:joint-excess}.
Equation \eqref{eq:eq-partner} gives
item~\ref{item:symmetric-equal-rank-partner}.
Finally, the first branch of \eqref{eq:one-sided-piecewise} gives
item~\ref{item:symmetric-small-comon-partner}.
\end{proof}

\begin{proof}[Proof of \cref{cor:symmetric-rank-nine-ten}]
Interchange the summands so that
$s(\tens{S}_1)\le9$ and $s(\tens{S}_2)\le10$.
Equations \eqref{eq:rank-nine-ten-excess} and
\eqref{eq:rank-ten-comon-one} give
$e_s(\tens{S}_1)\le4$ and $\Delta_c(\tens{S}_2)\le1$.
The first branch of \cref{prop:one-sided} therefore gives
$\Delta_s(\tens{S}_1,\tens{S}_2)=0$.
\end{proof}

\subsection{CP rank excess control of the symmetric rank gap}
\label{sec:cp-excess-symmetric-proof}

We now prove \cref{thm:sym-p}.  The case of small CP rank excess follows from
rank equality.  For the remaining case, we remove the cubes already
present in a minimum CP decomposition and bound the residual
symmetric tensor using its wedge supports.

The decrease of CP rank excess under deletion of displayed cubes
is proved over $\RR$ in \cite[Lemma~2.4]{li2026comonrankgap}.
The next proposition also controls the symmetric rank gap with
an arbitrary partner.

\begin{proposition}[Removing displayed cubes]
\label{prop:core-transfer}
Suppose that a minimum CP decomposition of a symmetric tensor
$\tens{S}_1$ consists of a symmetric subtotal $\tens{S}'_1$ and a list of
cubes.  Then
\[
0\le e_r(\tens{S}'_1)\le e_r(\tens{S}_1)
\]
and the following transfer inequality holds for every symmetric
partner $\tens{S}_2$:
\begin{equation}\label{eq:core-transfer}
\Delta_s(\tens{S}_1,\tens{S}_2)
\le
2\bigl(e_r(\tens{S}_1)-e_r(\tens{S}'_1)\bigr)
+\Delta_s(\tens{S}'_1,\tens{S}_2).
\end{equation}
\end{proposition}
\begin{proof}
By \cref{lem:minimum}, sublist minimality and support spanning give
\begin{align*}
r(\tens{S}_1)-r(\tens{S}'_1) 
=d(\tens{S}_1)+e_r(\tens{S}_1)-d(\tens{S}'_1)-e_r(\tens{S}'_1),\ \
r(\tens{S}_1)-r(\tens{S}'_1)\ge d(\tens{S}_1)-d(\tens{S}'_1),
\end{align*}
and hence $0\le e_r(\tens{S}'_1)\le e_r(\tens{S}_1)$.
Choose $d(\tens{S}_1)-d(\tens{S}'_1)$ removed cube directions whose images
form a basis modulo $\Ess(\tens{S}'_1)$.  Let $\tens{S}_{1,\mathrm{ind}}$ be their subtotal
and let $\tens{S}_{1,\mathrm{rem}}$ be the subtotal of the remaining cubes.
The chosen directions span a complement to $\Ess(\tens{S}'_1)$, so
\[
\tens{S}_1=\widetilde{\tens{S}}_1+\tens{S}_{1,\mathrm{rem}},\ \ 
\widetilde{\tens{S}}_1\eqdef\tens{S}'_1\oplus\tens{S}_{1,\mathrm{ind}},\ \ 
s(\tens{S}_{1,\mathrm{rem}})\le e_r(\tens{S}_1)-e_r(\tens{S}'_1).
\]
By \cref{lem:independent-cube},
the gaps satisfy
$\Delta_s(\widetilde{\tens{S}}_1,\tens{S}_2)=\Delta_s(\tens{S}'_1,\tens{S}_2)$.
Symmetric rank subadditivity also gives
\begin{align*}
s(\tens{S}_1)\le s(\widetilde{\tens{S}}_1)+e_r(\tens{S}_1)-e_r(\tens{S}'_1),\ \
s(\tens{S}_1\oplus\tens{S}_2)
\ge s(\widetilde{\tens{S}}_1\oplus\tens{S}_2)
-\bigl(e_r(\tens{S}_1)-e_r(\tens{S}'_1)\bigr).
\end{align*}
Subtracting proves \eqref{eq:core-transfer}.
\end{proof}

When the residual decomposition is symmetric in two modes,
the wedge budget gives a shorter estimate. The residual bounds
$d(\tens{S}'_1)\le2e_r(\tens{S}'_1)$ and
$r(\tens{S}'_1)\le3e_r(\tens{S}'_1)$ are the real estimates in
\cite[Lemma~4.6]{li2026comonrankgap}.

\begin{proposition}[A bound for partially symmetric CP decompositions]
\label{prop:partial-bound}
Suppose that $\tens{S}_1$ has a minimum CP decomposition of the form
$\sum_i\vect{a}_i\otimes\vect{v}_i^{\otimes2}$.  Then the following bound holds
for every symmetric partner $\tens{S}_2$:
\begin{equation*}
\Delta_s(\tens{S}_1,\tens{S}_2)
\le
\begin{cases}
2e_r(\tens{S}_1),&e_r(\tens{S}_1)\le3,\\[2pt]
9e_r(\tens{S}_1)-1,&e_r(\tens{S}_1)\ge4.
\end{cases}
\end{equation*}
\end{proposition}
\begin{proof}
If either tensor is zero, the gap is zero, so assume both are nonzero.
If $e_r(\tens{S}_1)\le3$, then \cref{prop:e3} gives
$e_s(\tens{S}_1)=e_r(\tens{S}_1)$, and \eqref{eq:elementary-excess}
proves the first branch.  Assume that $e_r(\tens{S}_1)\ge4$.
Remove every displayed cube and let $\tens{S}'_1$ be the residual.
By \cref{prop:core-transfer}, we have
$e_r(\tens{S}'_1)\le e_r(\tens{S}_1)$.
If $\tens{S}'_1=0$, that proposition gives a gap at most
$2e_r(\tens{S}_1)\le9e_r(\tens{S}_1)-1$. 
Otherwise, every remaining wedge $\vect{a}_i\wedge\vect{v}_i$ is
nonzero and its support contains both factors.
Equation \eqref{eq:wedge-budget} gives $d(\tens{S}'_1)\le2e_r(\tens{S}'_1)$ and
$r(\tens{S}'_1)\le3e_r(\tens{S}'_1)$.
The polarization identity for a product with a repeated factor gives
$s(\tens{S}'_1)\le3r(\tens{S}'_1)\le9e_r(\tens{S}'_1)$.
Combining \eqref{eq:minimum-rank-bound} with
\cref{prop:core-transfer} yields
\[
\Delta_s(\tens{S}_1,\tens{S}_2)
\le2\bigl(e_r(\tens{S}_1)-e_r(\tens{S}'_1)\bigr)+9e_r(\tens{S}'_1)-1
\le9e_r(\tens{S}_1)-1.
\]
The proof is complete. 
\end{proof}

We next combine the three estimates for wedge supports for a general
minimum CP decomposition.

\begin{proposition}[High CP rank excess bound]
\label{prop:12p}
If $e_r(\tens{S}_1)\ge4$, the following bound holds for every symmetric
partner $\tens{S}_2$:
\begin{equation*}
\Delta_s(\tens{S}_1,\tens{S}_2)\le12e_r(\tens{S}_1)-9.
\end{equation*}
\end{proposition}
\begin{proof}
The case $\tens{S}_2=0$ is immediate, so assume that the partner is nonzero.
Among all minimum CP decompositions of $\tens{S}_1$, choose one for which
$(\kappa_{ab},\kappa_{ca},\kappa_{bc})$ is lexicographically minimal,
and remove its displayed cubes.  Let $\tens{S}'_1$ be the residual
symmetric tensor and set $k\eqdef e_r(\tens{S}'_1)\le e_r(\tens{S}_1)$.
By \cref{lem:minimum}, the residual decomposition is minimum,
and deleting displayed cubes does not change any wedge span.
Its wedge-dimension triple is therefore lexicographically minimal among all minimum CP decompositions of $\tens{S}'_1$.
Indeed, adjoining the removed cubes to a minimum CP decomposition with a lexicographically smaller triple would contradict the original choice. 
If $k\le3$, then \cref{prop:e3} gives $e_s(\tens{S}'_1)=k$.
Equation \eqref{eq:elementary-excess} gives
$\Delta_s(\tens{S}'_1,\tens{S}_2)\le2k$ when $\tens{S}'_1\ne0$; the same
inequality is immediate when $\tens{S}'_1=0$.
Thus \cref{prop:core-transfer} gives
$\Delta_s(\tens{S}_1,\tens{S}_2)\le2e_r(\tens{S}_1)$.
If $k\ge4$ and $\tens{S}'_1$ has a minimum CP decomposition symmetric
in two modes, then \cref{prop:partial-bound,prop:core-transfer}
give $$\Delta_s(\tens{S}_1,\tens{S}_2)\le9e_r(\tens{S}_1)-1.$$
Both bounds are at most $12e_r(\tens{S}_1)-9$ for $e_r(\tens{S}_1)\ge4$. 

It remains to consider $k\ge4$ when no minimum partially symmetric CP decomposition
exists.  For each mode pair, at least one connected wedge component has
matroid rank at least three.  Otherwise all components could be replaced
as in the proof of \cref{lem:small-components}, producing a
minimum partially symmetric CP decomposition.  \Cref{lem:connected-wedge} therefore gives
\begin{equation}\label{eq:nonpartial-two-saving}
h_{xy}\le2\kappa_{xy}-2
\ \ (xy\in\{ab,ca,bc\}).
\end{equation}

Each wedge dimension is at most $k$ by \eqref{eq:wedge-budget}.
If $\kappa_{ab}=k$,
then every minimum CP decomposition of $\tens{S}'_1$ has $ab$-wedge dimension $k$;
the chosen decomposition therefore minimizes $\kappa_{ca}$ without a
side condition.  If also $\kappa_{ca}=k$, it similarly minimizes the
third wedge dimension.  Apply \cref{lem:small-components} to each mode pair whose wedge
dimension is minimized unconditionally, and use
\eqref{eq:nonpartial-two-saving} for the other pairs.
Put $h\eqdef h_{ab}+h_{ca}+h_{bc}$.
If $\kappa_{ab}<k$, these estimates give
$h\le4(k-1)/3+4k-4$.
If $\kappa_{ab}=k$ and $\kappa_{ca}<k$, they give
$$h\le4k/3+4(k-1)/3+2k-2.$$ 
If both dimensions equal $k$, all three minima are unconditional
and $h\le4k$.
Each of these bounds is at most $16(k-1)/3$ for $k\ge4$.
By \cref{lem:three-sector}, we have $d(\tens{S}'_1)\le h/2$, and hence
\[
s(\tens{S}'_1)-1\le4k+\frac32h-1\le12k-9.
\]
Finally, \cref{prop:core-transfer} and
\eqref{eq:minimum-rank-bound} yield
\begin{align*}
\Delta_s(\tens{S}_1,\tens{S}_2)
\le2\bigl(e_r(\tens{S}_1)-k\bigr)+s(\tens{S}'_1)-1 
\le2\bigl(e_r(\tens{S}_1)-k\bigr)+12k-9 
\le12e_r(\tens{S}_1)-9.
\end{align*}
\end{proof}

\begin{proof}[Proof of \cref{thm:sym-p}]
Under \eqref{eq:smaller-cp-excess}, suppose first that $e_r(\tens{S}_1)\le3$.
Then \cref{prop:e3} gives $r(\tens{S}_1)=s(\tens{S}_1)$ and
$e_s(\tens{S}_1)=e_r(\tens{S}_1)$.
The first branch of \eqref{eq:symmetric-gap-from-cp-excess}
therefore follows from \eqref{eq:elementary-excess}.
For $e_r(\tens{S}_1)\ge4$, apply \cref{prop:12p}. 
Finally, if $e_r(\tens{S}_2)\le3$, then
\eqref{eq:smaller-cp-excess} gives $e_r(\tens{S}_1)\le3$ as well.
Both summands therefore have equal CP and symmetric ranks.
Their CP rank is additive by
\eqref{eq:ordinary-small-excess-additivity}.  Hence
\[
s(\tens{S}_1)+s(\tens{S}_2)
=r(\tens{S}_1\oplus\tens{S}_2)
\le s(\tens{S}_1\oplus\tens{S}_2)
\le s(\tens{S}_1)+s(\tens{S}_2),
\]
which proves \eqref{eq:cp-excess-three-symmetric-additivity}.
\end{proof}

% \begin{remark}[Relation to CP rank criteria]
% Friedland's criterion for equality of CP rank and symmetric rank
% \cite[Theorem~1.1]{friedland2016remarks} and the CP rank additivity
% criterion in \cite[Theorem~1.1]{buczynski2020additivity}
% give symmetric rank additivity in a smaller range of CP rank excess.
% The estimates above also bound a positive symmetric rank gap.
% Their additional step is to remove displayed cubes and control the
% remaining symmetric tensor through the wedge supports of a minimum
% CP decomposition.
% \end{remark}

\subsection{The shifted maximum symmetric rank bound}\label{sec:shifted-proof}
We now turn to essential dimension.  The following consequence of the
Alexander--Hirschowitz theorem completes the polar extraction argument
from \cref{lem:polar-leading}.

\begin{proposition}[Regular cube completion]
\label{lem:regular-cube-completion}
Let $\dim\set{H}=m\ge1$. Then $g_3(m)\ge m$.
Every ordered linearly independent family
$\vect{\ell}_1,\ldots,\vect{\ell}_q\in\set{H}$ with $0\le q\le m$
extends to a tuple
$\vect{\ell}_1,\ldots,\vect{\ell}_{g_3(m)}$ satisfying
\[
\rank D\Phi_{g_3(m)}
(\vect{\ell}_1,\ldots,\vect{\ell}_{g_3(m)})
=\dim\Sym^3(\set{H})=\binom{m+2}{3}.
\]
\end{proposition}

\begin{proof}
The real assertion is \cite[Lemma~6.8]{li2026comonrankgap}.
Put $g\eqdef g_3(m)$. Formula \eqref{eq:generic} gives $g\ge m$.
The Alexander--Hirschowitz theorem
\cite[Theorem~2]{bernardi2018hitchhiker} implies that the complex
map $\Phi_g$ is dominant.
In characteristic zero, its differential is generically surjective.
Thus its Jacobian has a nonzero maximal minor $J$ with rational
coefficients. 
Let $K$ be a nonzero $q\times q$ coordinate minor of the first
$q$ vector variables, with $K=1$ when $q=0$.
Since $JK$ is a nonzero polynomial and $\FF$ is infinite,
there exists
$\vect{u}=(\vect{u}_1,\ldots,\vect{u}_g)\in\set{H}^g$
such that $J(\vect{u})K(\vect{u})\ne0$.
Hence $D\Phi_g(\vect{u})$ is surjective and
$\vect{u}_1,\ldots,\vect{u}_q$ are linearly independent. 
Choose $L\in\mathrm{GL}(\set{H})$ such that
$L\vect{u}_i=\vect{\ell}_i$ for $1\le i\le q$.
Differentiating the equivariance identity for $\Phi_g$ gives
\[
D\Phi_g(L^{\times g}\vect{u})\circ L^{\times g}
=\Sym^3(L)\circ D\Phi_g(\vect{u}).
\]
Both $L^{\times g}$ and $\Sym^3(L)$ are invertible.
Therefore $D\Phi_g(L^{\times g}\vect{u})$ is surjective,
and the tuple $L^{\times g}\vect{u}$ gives the required completion.
\end{proof}

\begin{proof}[Proof of \cref{prop:shifted}]
Let $F\in\Sym^3(\set{V})$ be nonzero with $\dim\set{V}=d\ge2$.
Choose $\alpha\in\set{V}^*$ such that $F(\alpha)\ne0$ and set
$\set{H}\eqdef\ker\alpha$,
$h\eqdef\rank C_F(\alpha).$ 
By \cref{lem:polar-leading}, there are sums $S_\varepsilon$ of $h$ cubes such that
$R_\varepsilon\eqdef F-S_\varepsilon\in\Sym^3(\set{H})$ and
$\varepsilon R_\varepsilon\to-P$, where $P$ is a sum of cubes of
$h-1$ linearly independent linear forms in $\set{H}$. 
Complete those $h-1$ cube factors of $-P$ to a regular $g_3(d-1)$-tuple by
\cref{lem:regular-cube-completion}. Append the negatives of its
last $g_3(d-1)-h+1$ entries. 
Regularity follows because the differential with respect to the first $g_3(d-1)$ entries is already surjective. 
Apply the regular neighborhood assertion of
\cref{lem:regular-degeneration} in $\Sym^3(\set{H})$.
For sufficiently small nonzero $\varepsilon$, the cubic
$\varepsilon R_\varepsilon$ therefore has symmetric rank at most
$2g_3(d-1)-h+1$.  Nonzero scaling preserves symmetric rank.
Restoring the $h$ extracted cubes gives
\[
s(F)\le h+s(R_\varepsilon)
\le h+2g_3(d-1)-h+1=2g_3(d-1)+1.
\]
For $1\le d\le4$, use \eqref{eq:rank-cap}. 
Finally, \cref{lem:elementary-one-sided} gives
$$\Delta_s(\tens{S}_1,\tens{S}_2)\le s(\tens{S}_1)-1
\le M_{3,\FF}(d(\tens{S}_1))-1.$$
\end{proof}

% \begin{remark}[Relation to maximum symmetric rank bounds]
% The Alexander--Hirschowitz theorem determines $g_3(n)$; see
% \cite[Theorem~2]{bernardi2018hitchhiker}.  Blekherman and Teitler
% prove the general estimate $M_{3,\FF}(n)\le2g_3(n)$ over both fields
% \cite[Theorem~1 for $\CC$ and Theorems~2--3 for $\RR$]
% {blekherman2015maximum}, while Jelisiejew proves
% $M_{3,\CC}(n)\le\binom{n+1}{2}-1$ for $n\ge3$
% \cite[Corollary~6, p.~331]{jelisiejew2014upper}.  The real shifted bound is
% \cite[Proposition~6.6]{li2026comonrankgap}.
% The proof above uses the same polar deformation and regular completion
% over both fields. This yields the dimension bound in
% \cref{thm:sym-d}.
% \end{remark}

\begin{proof}[Proof of \cref{thm:sym-d}]
By \cref{prop:shifted}, the gap is at most
$M_{3,\FF}(d(\tens{S}_1))-1$.
The part concerning small dimensions of that proposition gives
the first branch of \eqref{eq:sym-d-main}; substituting
\eqref{eq:generic} into its shifted bound gives the other two branches,
with the exceptional value at $d(\tens{S}_1)=6$.
The case $d(\tens{S}_1)=1$ gives
\eqref{eq:symmetric-dimension-one-additivity}.
If $d(\tens{S}_1),d(\tens{S}_2)\le4$, then \eqref{eq:rank-cap} gives
\[
e_s(\tens{S}_1)\le d(\tens{S}_1)-1\le3,
\ \ 
e_s(\tens{S}_2)\le d(\tens{S}_2)-1\le3.
\]
Thus item~\ref{item:symmetric-excess-rectangle} of \cref{thm:sym-q}
gives \eqref{eq:symmetric-dimension-four-additivity}.
\end{proof}

\begin{remark}
For $\FF=\RR$, \cite[Corollary~6.11]{li2026comonrankgap} gives
$M_{3,\RR}(6)\le16$.
Together with \cref{prop:shifted}, this implies
$\Delta_s(\tens{S}_1,\tens{S}_2)\le15$ whenever
$\min\{d(\tens{S}_1),d(\tens{S}_2)\}=6$.
\end{remark} 

\begin{proof}[Proof of \cref{cor:symmetric-dimension-asymptotics}]
For $d(\tens{S}_1)\ge5$ and $d(\tens{S}_1)\ne6$, one has
\[
2\left\lceil\frac{d(\tens{S}_1)(d(\tens{S}_1)+1)}6\right\rceil
\le\frac{d(\tens{S}_1)(d(\tens{S}_1)+1)}3+2.
\]
When $d(\tens{S}_1)=6$, the exceptional bound equals
$16=d(\tens{S}_1)(d(\tens{S}_1)+1)/3+2$, and the bounds for
$1\le d(\tens{S}_1)\le4$ are smaller.
This proves \eqref{eq:symmetric-uniform-dimension}.
\end{proof}

\section{Conclusions}
\label{sec:conclusion}

We have established quantitative bounds for the Strassen CP rank gap
and symmetric rank gap of third-order tensors over $\RR$ and
$\CC$.  The bounds depend on CP rank excess, symmetric rank excess,
the corresponding ranks, and maximum mode ranks (essential dimensions in the symmetric case).  They also give
sufficient conditions for CP rank or symmetric rank additivity under
direct sums.
This naturally leads to the following question.

\begin{problem}
Can the control functions in the gap bounds of this paper be
sharpened over $\RR$ and $\CC$? In particular, can the quadratic
dependence on maximum mode rank (essential dimension in the symmetric case) be reduced to subquadratic or linear dependence,
and can the bounds in terms of CP rank excess and symmetric rank
excess be improved?
\end{problem}

\noindent\textbf{Acknowledgments.} 
The author acknowledges the assistance of ChatGPT 6 Astra throughout the mathematical development and preparation of this manuscript. 
The author takes full responsibility for the content of this paper. 
This work was supported by the Guangdong Basic and Applied Basic
Research Foundation (2026A1515050001, 2026A1515010188).

\bibliographystyle{siamplain}
\bibliography{ReferenceTensor}
\end{document}